\documentclass[11pt]{amsart}
\usepackage{amsmath}
\usepackage{amssymb}
\usepackage{amsthm}
\usepackage{bbm}
\usepackage{graphicx}
\usepackage{algorithm}
\usepackage{algorithmic}
\usepackage{clrscode}
\usepackage{fullpage}
\usepackage{stmaryrd}
\usepackage{color}
\usepackage{url}
\usepackage{caption}
\usepackage{subcaption}

\newtheorem{cor}{Corollary}
\newtheorem{thm}{Theorem}
\newtheorem{lem}{Lemma}

\newtheorem{Def}{Definition}

\DeclareMathOperator*{\argmin}{arg\,min}

\begin{document}
\pagestyle{plain}
\newenvironment{frcseries}{\fontfamily{frc} \selectfont}{}
\newcommand{\textfrc}[1]{{\frcseries #1}}
\newcommand{\mathfrc}[1]{\text{\textfrc{#1}}}

\title{Rapidly Computing Sparse Legendre Expansions via Sparse Fourier Transforms}\thanks{Xianfeng (Janice) Hu:  Institute for Mathematics and its Applications, University of Minnesota ({\tt xhu@umn.edu}).\\
\indent M.A. Iwen:  Department of Mathematics and Department of ECE, Michigan State University ({\tt markiwen@math.msu.edu}).  M.A. Iwen was supported in part by NSF DMS-1416752.\\ 
\indent Hyejin Kim:  Department of Mathematics and Statistics, University of Michigan -- Dearborn ({\tt khyejin@umich.edu}).}
%
\author{
{Xianfeng Hu}
\and
{Mark Iwen} %
\and
{Hyejin Kim}
}

\begin{abstract}
In this paper we propose a general strategy for rapidly computing sparse Legendre expansions.  
The resulting methods yield a new class of fast algorithms capable of  approximating a given function $f: [-1,1] \rightarrow \mathbbm{R}$ with a near-optimal linear combination of $s$ Legendre polynomials of degree $\leq N$ in just $(s \log N)^{\mathcal{O}(1)}$-time.  When $s \ll N$ these algorithms exhibit sublinear runtime complexities in $N$, as opposed to traditional $\Omega(N \log N)$-time methods for computing all of the first $N$ Legendre coefficients of $f$.  Theoretical as well as numerical results demonstrate the effectiveness of the proposed methods. 
\end{abstract}

\maketitle
\thispagestyle{empty}

\section{Introduction}
\label{sec:Intro}

In this paper we consider {\it Legendre-compressible functions} which can be well approximated by a linear combination of a small number of unknown, and potentially high-degree, Legendre polynomials.  
Given such a function our objective is to quickly learn the best basis of Legendre polynomials with which to approximate it, and then to compute their coefficients.  Let $f: [-1,1] \rightarrow \mathbbm{R}$ be a degree $N$ polynomial, and $L_n(x)$ denote the Legendre polynomial of degree $n$.  We aim to rapidly and accurrately compute $f$'s Legendre coefficients, $\tilde{f}(n) \in \mathbbm{R}$ for $n \in \{0, \dots, N\}$ with
\begin{equation}
f(x) = \sum^N_{n=0} \tilde{f}(n) L_n(x),
\label{Probf}
\end{equation}
whenever $\tilde{f}(n) \approx 0$ for all but $s \ll N$ {\it initially unknown} values of $n$.  We will call any numerical method with this objective a {\it sparse Legendre expansion algorithm}.  

Note that solving this problem is straightforward if one is willing to sample $f$ at $N+1$ points in $[-1,1]$, and then compute all $N+1$ of its Legendre coefficients.  However, any such approach will necessarily require $\Omega(N)$-operations, which can become overwhelming when the maximal degree, $N$, of $f$ is large. Our objective here is to select the best basis of $s \ll N$ Legendre polynomials of degree $\leq N$ for $f$, and then estimate their coefficients, in $(s \log N)^{\mathcal{O}(1)}$-time.  When $s$ is significantly smaller than $N$, these methods will be faster than any traditional approach which computes all $N$ Legendre coefficients of $f$.  Fast sparse Legendre expansion algorithms of this kind are a natural first step toward the development of computationally tractable algorithms for approximating functions of many variables with respect to tensorized Legendre polynomial bases.  In such multivariate problems the maximal degree, $N$, grows exponentially in the number of variables, rapidly rendering even $\mathcal{O}(N)$-time methods impractical.  Our longterm goal is to extend $(s \log N)^{\mathcal{O}(1)}$-time sparse expansion methods for functions one variable, once they are properly understood, to the multi-variate setting.  If possible, the resulting methods would be of value in many computational applications including, e.g., uncertainty quantification \cite{LeMatre2010spectral} and the computation of polynomial chaos expansions \cite{cohen2010convergence,cohen2011analytic}. 

The majority of previously proposed sparse Legendre expansion methods are based on Prony-like approaches.  Examples include results by Peter et. al. \cite{peter2013representation} who develop a method based on a more general approach from \cite{peter2013generalized} which needs only $\mathcal{O}(s)$ samples from (various derivatives of) $f$ in order to recover its $s$-sparse Legendre expansion.  More recently, Potts and Tasche \cite{pottsLegendre} used approximation techniques to adapt previous Prony-like methods for the recovery of Chebyshev-sparse functions \cite{potts2014sparse} to the Legendre-sparse setting.  However, no theoretical results are proven in \cite{pottsLegendre} that demonstrate the methods therein can extend to functions with compressible (as opposed to exactly $s$-sparse) Legendre expansions, nor is it proven that they can tolerate even modest levels of noise in general.  
In contrast, herein we provide a theoretical support recovery guarantee which proves that our techniques can indeed locate the principle support of a relatively large class of Legendre-compressible functions (see, e.g., Theorem~\ref{thm:SupportIdentification}).  When combined with coefficient estimation methods based on techniques from compressive sensing (see, e.g., Lemma~\ref{lem:CoefEstLemma}) these support recovery guarantees allow one to prove a variety of general sublinear-time recovery guarantees for functions with compressible Legendre expansions.

Other sparse Legendre expansion methods include those based on compressive sensing approaches \cite{HolgerBook}.  In particular, Rauhut and Ward \cite{rauhut2012sparse} demonstrate that $\mathcal{O}(s \cdot \log^4 N)$ samples from $f$ suffice in order to accurately and stably approximate $f$ with a near-optimal sparse Legendre expansion.  The associate reconstruction algorithms are $\Omega(N)$-time, however.  As opposed to these previous sparse Legendre expansion methods, we propose a new approach motivated by a recently proposed FFT-based algorithm for computing all $N+1$ Legendre coefficients of a given function $f$ in $\mathcal{O}(N \log N)$-time.  In the process, we demonstrate a general approach which allows one to utilize any Sparse Fourier Transform algorithm (see, e.g., \cite{GMS2005Improved,HIKP2012Simple,LWC2013adaptive,SI2013Improved}) one desires in order to recover functions which are sparse/compressible in other polynomial bases (herein, Chebyschev and Legendre).  

Computing all $N+1$ Legendre coefficients of $f$ in $o(N^2)$-time is itself a challenging problem which has attracted a good deal of attention over the past two decades.  Proposed methods include, e.g., fast multipole-like approaches \cite{alpert1991fast}, 
and algorithms based on integral transform techniques \cite{de2013expansion}, to mention just a few.  These methods are $\mathcal{O}(N \log^c N)$-time for various $c \in \mathbbm{R}^+$.  Most pertinent to the sparse Legendre expansion algorithms proposed herein, however, are recent FFT-based algorithms for computing all $N$ Legendre coefficients in $\mathcal{O}(N \log N)$-time \cite{iserles2011fast}.  These methods work by $(i)$ implicitly mapping the Legendre coefficients of $f$ to the Fourier coefficients of a related function, $f_r$, which remains easily to sample, followed by $(ii)$ computing the Fourier coefficients of $f_r$ with an FFT, and then $(iii)$ using the computed Fourier coefficients of $f_r$ in order to recover the Legendre coefficients of the original function $f$ (i.e., by inverting the map from $(i)$).  The sparse Legendre expansion methods proposed herein are based on this same type of approach.  In particular, the algorithms proposed herein result from combining ideas from Iserles' FFT-based Legendre algorithm \cite{iserles2011fast} (summarized below in \S \ref{sec:IserlesMap}) with sparse Fourier transform techniques (briefly discussed in the next section).  

\subsection{Sparse Fourier Transforms}
\label{sec:SFT}

Sparse Fourier Transforms (SFTs) are algorithms for quickly computing near-optimal sparse approximations to the Fourier series of a given periodic function $f:[-\pi, \pi]^D \rightarrow \mathbb{C}$.  Suppose that $f$ is a trigonometric polynomial of degree $N$ in every variable so that the Fourier series of $f$ is effectively $\hat{f} \in \mathbb{C}^{N^D}$.  An optimal $s$-term trigonometric approximation to $f$ is given by
\begin{equation}
f^{\rm opt}_s ({\bf x}) := \sum^s_{j = 1} \hat{f} \left( {\boldsymbol \omega}_j \right) \mathbbm{e}^{\mathbbm{i} {\boldsymbol \omega}_j \cdot {\bf x} }
\label{BestkTerm}
\end{equation}
where ${\boldsymbol \omega}_1, \dots, {\boldsymbol \omega}_{N^D} \in (-N/2, N/2]^D \cap \mathbb{Z}^D$ are ordered by the magnitudes of their Fourier coefficients, $\hat{f} \left( {\boldsymbol \omega}_j \right) \in \mathbb{C}$, so that
\begin{equation}
\big| \hat{f}({\boldsymbol \omega}_1) \big| \geq \big| \hat{f}( {\boldsymbol \omega}_2 ) \big| \geq \dots \geq \big| \hat{f}({\boldsymbol \omega}_{N^D}) \big|.
\label{SortedCoefs}
\end{equation}
The optimal $s$-term approximation error is then $\| f - f^{\rm opt}_s \|_2 = \| \hat{f} - \hat{f}^{\rm opt}_s \|_2$.  In this setting, any discrete Fourier method will take function evaluations of $f$ as input, and then output an approximate $f^{\rm opt}_s$ for some value of $s \in \{ 0, \dots, N^D \}$.  A standard FFT always uses $s = N^D$, and so recovers trigonometric polynomials exactly.  Sparse FFTs allow $s$ to be chosen independently of $N^D$, and exactly recover all trigonometric polynomials consisting of at most $s$ nonzero terms.

The primary objective of SFTs is to compute an accurate approximation to $\hat{f}$, and therefore to $f$ itself, as quickly as absolutely possible.  This has lead to the development of a wide range of algorithms for approximating degree $N$ trigonometric polynomials which use $(s \cdot \log N)^{\mathcal{O}(1)}$ floating point operations, where $s$ is the user specified sparsity parameter.  In order to achieve these operation counts, all such SFTs can utilize at most $(s \cdot \log N)^{\mathcal{O}(1)}$ function evaluations from $f$ during their execution.  It is interesting to note, for the purposes of comparison, that SFTs are therefore closely related to Fourier-based compressed sensing techniques \cite{CT2005Decoding,D2006Compressed,CT2006Near,CRT2006Stable,RV2006Sparse,HolgerBook} whose primary objective is to approximate $f$ using as few function evaluations, or samples, as absolutely possible (see, e.g., \cite{BIS2012Design}).  However, although the relationship between SFTs and compressed sensing is very close, SFTs generally utilize more samples than the best compressed sensing methods in practice.  Similarly, even the fastest Fourier-based compressed sensing methods are too slow to serve as SFTs since SFTs purposefully exceed the strictest sampling requirements of compressed sensing methods in order to reduce their runtime complexities as much as possible.

The first sparse Fourier methods were essentially approximate Hadamard transforms that were developed by researchers in the machine learning community for quickly learning boolean functions of many variables (see, e.g., \cite{KM1993Learning, BFJKMR1994Weakly,M1994Learning} and \cite{GL1989Hard,G1999Foundations}).  These techniques were later adapted to produce randomized SFTs for approximating trigonometric polynomials as rapidly as possible \cite{M1995Randomized,GGIMS2002Near,AGS2003Proving,GMS2005Improved}.   These subsequent SFTs all take random samples of a given periodic function $f$ as input, and then output a trigonometric polynomial, $y$, of degree $N$ which satisfies $\| f - y \|_2 \approx \| f - f^{\rm opt}_s \|_2$ with high probability.  The fastest of these SFTs~\cite{GMS2005Improved} uses only $s \cdot \log^{\mathcal{O}(1)} N$ operations.  As a result, it is faster than the FFT for accurately approximating periodic functions which are dominated by $s \ll N$ of their largest magnitude Fourier coefficients \cite{AAFFT1exp}.

Over the last several years SFTs have been improved significantly in both theory and practice.  Recent work includes better implementations \cite{HIKP2012Simple,LWC2013adaptive,SI2013Improved}, improvements in runtime complexity bounds (both upper and lower) \cite{HIKP2012Nearly,LWC2013adaptive,I2013Improved}, adaptation of the methods to the recovery of superpositions of sinusoids with non-integer frequencies \cite{BCGLS2012What}, and improvements in theoretical error guarantees \cite{A2010Deterministic,I2010Combinatorial,SI2013Improved,I2013Improved}.  In particular, entirely deterministic SFTs exist \cite{I2010Combinatorial,I2013Improved} that are guaranteed to \textit{always} return a near-optimal sparse trigonometric polynomial, $y_s: [-\pi, \pi]^D \rightarrow \mathbb{C}$, having
$\| f - y_s \|_2 \approx \| f - f^{\rm opt}_s \|_2$.  More specifically, the following theorem was proven in \cite{I2013Improved}.

\begin{thm}
Suppose $f: [-\pi,\pi]^D \rightarrow \mathbb{C}$ has $\hat{f}(\omega_1, \dots, \omega_D) = 0$ if $\left( \omega_1, \dots, \omega_D \right) \notin$ $\left( [-\frac{N}{2}, \frac{N}{2}] \cap \mathbb{Z} \right)^D$.
Let $s, \epsilon^{-1} \in \mathbb{N} \setminus \{ 1 \}$ with $(s/\epsilon)^2 \geq 4$.  Then, there exists a simple deterministic algorithm that is guaranteed to output a trigonometric polynomial, $y_s: [-\pi, \pi]^D \rightarrow \mathbb{C}$,
satisfying
\begin{equation}
\left\| f - y_s \right\|_2 \leq \left\| \hat{f} - \hat{f}^{\rm~ opt}_s \right\|_2 + \frac{22\epsilon \cdot \left\| \hat{f} - \hat{f}^{\rm~ opt}_{(s/\epsilon)} \right\|_1}{\sqrt{s}}.
\label{eqn:MultReconError}
\end{equation}
The algorithm's operation count is
\begin{equation}
\mathcal{O} \left( \frac{s^2 \cdot D^4 \cdot \log^4 (ND)}{\log \left(\frac{s}{\epsilon} \right) \cdot \epsilon^2} \right).
\label{eqn:DetRuntime}
\end{equation}

If succeeding with probability $(1-\delta) \in [2/3,1)$ is sufficient, and $(s/\epsilon) \geq 2$, a Monte Carlo variant of the deterministic algorithm may be used.  This Monte Carlo variant will output a trigonometric
polynomial, $y_s: [-\pi, \pi]^D \rightarrow \mathbb{C}$, that satisfies Equation~\ref{eqn:MultReconError} with probability at least $1-\delta$.  Its operation count will be
\begin{equation}
\mathcal{O} \left( \frac{s \cdot D^4}{\epsilon} \cdot \log^3 (ND) \cdot \log \left( \frac{ND}{\delta} \right) \right).
\label{eqn:RandRuntime}
\end{equation}
\label{thm:MultDimRecov}
\end{thm}

The Fourier algorithms referred to by Theorem~\ref{thm:MultDimRecov} are able to accurately approximate the discrete Fourier transform of a given function much more quickly than standard Fast Fourier Transform (FFT) methods \cite{B2001Chebyshev} whenever the sorted magnitudes of the Fourier coefficients \eqref{SortedCoefs} go to zero quickly enough \cite{SI2013Improved}.  More specifically, the developed Fourier approximation algorithms have operation counts that scale \textit{polynomially} in $D$ and $\log N$, as opposed to standard FFT methods whose operation counts scale \textit{exponentially} in $D$ and $\log N$.  We direct the reader to \cite{gilbert2014recent} for a recent survey of SFT techniques, as well as for an easy introduction to their design and implementation.

\subsection{Our Proposed Sparse Legendre Expansion Method}
\label{sec:OurMethod}

The thought behind our proposed sparse Legendre expansion approach is naively simple at first glance.  One thinks:  {\it ``FFT-based algorithms for computing all $N$ Legendre coefficients like \cite{iserles2011fast} appear to work well.  Maybe SFTs can replace the FFTs in these methods in order to allow us to rapidly approximate Legendre-compressible functions!"}  Of course, a multitude of technical difficulties present themselves almost immediately.  Mainly, mapping the sparse Legendre coefficients of $f$ to a set of Fourier coefficients of a related function, $f_r$, is only helpful for computing sparse Legendre expansions if the map preserves sparsity.  The speed and accuracy of SFT methods depend on the sparsity of the function to which they are applied.  If the sparse Legendre coefficients of $f$ don't map to a set of Fourier coefficients of $f_r$ that are also fairly sparse, then SFT methods will not be able to approximate $f_r$ quickly enough to be interesting.  Furthermore, and perhaps more obviously, the map from the Legendre coefficients of $f$ to the Fourier coefficients of $f_r$ must be fairly well-conditioned in the $\ell_{\infty}$-sense.  If the map sends a few of the large-magnitude Legendre coefficients of $f$ to Fourier coefficients of comparatively tiny magnitude, then we will have difficultly recovering them with an SFT.  Finally, the inverse map from the computed Fourier coefficients of $f_r$ back to the Legendre coefficients of $f$ must be fast (e.g., well approximated by a sparse matrix multiply).  If not, we will not be able to use our computed Fourier coefficients of $f_r$ in order to compute the Legendre coefficients of the original function $f$ quickly enough to be of interest.  

Unfortunately, achieving all of these properties at once appears to be quite difficult.  Herein we take advantage of the fact that Iserles' map from Legendre coefficients to Fourier coefficients is fairly well-behaved with respect to sparsity (see \S \ref{sec:IserlesMap} and \S \ref{sec:ErrorAnal}).  In particular, we show that Legendre-sparse functions are mapped to Fourier-compressible functions (i.e., the map preserves sparsity fairly well).  Unfortunately, it appears to be impossible to force the map to also be both well-conditioned in the $\ell_{\infty}$-sense, {\it and} quickly invertible.  To compensate for this defect we modify our initial idea and employ a two-stage approach instead:  We first use a ``pretty well-conditioned'' version of Iserles' map in combination with SFT methods in order to rapidly identify the large-magnitude Legendre coefficients in $f$, and then use results form compressive sensing in order to accurately  approximate the identified Legendre coefficients.  Doing so allows us to develop workable sparse Legendre expansion algorithms which run in sublinear-time for sufficiently small sparsities.  We refer the reader to \S \ref{sec:Alg} for additional details, and to the next section for a related example.

\subsection{A Simple Example:  Recovering Sparse Chebyshev Expansions via SFTs}

In this section we briefly consider Chebyshev-sparse functions of the form
$$g(x) = \sum_{n \in \mathcal{S}} a_n T_n(x),~~\mathcal{S} \subset \{0, \dots, N \}, ~~|\mathcal{S}| = s \ll N,$$
where $T_n(x)$ denotes the degree $n$ Chebyshev polynomial of the first kind.  In this setting it is beneficial to consider the standard transformation $h = g( \cos x)$.  It is well known that this transformation implicitly maps the $n^{\rm th}$ Chebyshev coefficient of $g$ to the $n^{\rm th}$ Fourier cosine series coefficient of $h$ (see, e.g., \cite{B2001Chebyshev}).  In particular, we have that
$$h(x) = \sum_{n \in \mathcal{S}} a_n T_n(\cos x) = \sum_{n \in \mathcal{S}} a_n \cos (n x) = \sum_{n \in \mathcal{S}} \frac{a_n}{2} \left( \mathbbm{e}^{\mathbbm{i} n x} + \mathbbm{e}^{-\mathbbm{i} n x} \right).$$
It is now straightforward to see that creating $h$ by resampling $g$ according to the cosine function implicitly produces a sparsity-preserving linear map from the Chebyshev coefficients of $g$ to the positive Fourier coefficients of $h$,
\begin{equation}
\left( \begin{array}{l}  \hat{h}(0) \\ ~~\vdots \\  \hat{h}(N) \end{array} \right) = \left( \begin{array}{lllll}  \frac{1}{2} & 0 & \dots & 0 & 0\\ 
0 & \frac{1}{2} & 0 & \dots & 0 \\
0 & 0 & \ddots & 0 & 0\\
0  & \dots & 0 & \frac{1}{2} & 0\\
0 & 0 & \dots & 0 & \frac{1}{2}
\end{array} \right) \left( \begin{array}{l}  a_0 \\ ~~~~\vdots \\ a_N \end{array} \right).
\label{equ:ChebMap}
\end{equation}

Note that the map above \eqref{equ:ChebMap} has all of the desirable properties mentioned in section~\ref{sec:OurMethod}.  It is exactly sparsity preserving, well conditioned, and trivially invertible.  As a consequence, one can easily compute the sparse Chebyshev expansion of $g$ using SFTs.  One simply applies the SFT of their choice to $h$ and then reconstructs the Chebyshev coefficients of $g$ from the result using their knowledge of \eqref{equ:ChebMap}.  One goal of the research initiated here is to find a good sparsity-preserving map similar to \eqref{equ:ChebMap} for use with Legendre polynomials, if possible.  In this paper we make a first attempt toward this goal by analyzing the maps proposed by Iserles in \cite{iserles2011fast}.

\section{Notation and Background}
\label{sec:background}

We will denote the Legendre polynomial of degree $n$ by $L_n$.  The $n^{\rm th}$ Legendre coefficient of $f:  [-1,1] \rightarrow \mathbb{R}$ is then
\begin{equation}
\tilde{f}(n) := \left( n + \frac{1}{2} \right) \int^1_{-1} f(x) L_n(x)~dx
\end{equation}
for each $n \in \mathbb{N}$.  The Fourier series coefficients of a function $f:  [-\pi,\pi] \rightarrow \mathbb{C}$ will be denoted by
\begin{equation}
\hat{f}(\omega) := \frac{1}{2 \pi}\int^{\pi}_{-\pi} f(x) \mathbbm{e}^{-\mathbbm{i} \omega x}~dx
\end{equation}
for all $\omega \in \mathbb{Z}$.  The sequence of Legendre or Fourier coefficients of an appropriate function $f$ will be called $\tilde{f}$ or $\hat{f}$, respectively.

For any matrix $X \in \mathbb{C}^{N \times n}$ we will denote the $j^{\rm th}$ column of $X$ by ${\bf X}_j \in \mathbb{C}^{N}$.  The adjoint of a matrix, $X \in \mathbb{C}^{N \times n}$, will be denoted by $X^{*} \in \mathbb{C}^{n \times N}$, and the singular values of any matrix $X \in \mathbb{C}^{N \times n}$ will always be ordered  as $\sigma_1(X) \geq \sigma_2(X) \geq \dots \geq \sigma_{\min(N,n)}(X) \geq 0.$  Also, the condition number of the matrix $X$ will denoted by $\kappa(X) := \sigma_1(X) / \sigma_{\min(N,n)}(X)$.  We will use the notation $[N] := \{ 0, \dots, N\} \subset \mathbb{N}$ for any $N \in \mathbb{N}$.  For any matrix $X \in \mathbb{C}^{N \times (n+1)}$ and set $\mathcal{S} \subset [n]$ the matrix $X_{\mathcal{S}} \in \mathbb{C}^{N \times |\mathcal{S}|}$ will be the submatrix of $X$ formed by selecting the columns of $X$ indexed by $\mathcal{S}$.  Similarly, for any vector ${\bf x} \in \mathbbm{C}^N$ and set $\mathcal{S} \subset [n]$ the vector ${\bf x}_{S} \in \mathbbm{C}^N$ will have entries
$$(x_{S})_{j} =
\begin{cases}
0, & \text{if } j\notin S \\
x_j, & \text{if } j\in S\\
\end{cases}.$$
Finally, given any ${\bf x} \in \mathbb{C}^N$, the vector ${\bf x}^{\rm~ opt}_s \in \mathbb{C}^N$ will always denote an optimal $s$-sparse approximation to ${\bf x}$.  That is, ${\bf x}^{\rm~ opt}_s$ will always (i) have at most $s \in [N]$ nonzero entries, and (ii) satisfy
\begin{equation}
\| {\bf x} - {\bf x}^{\rm~ opt}_s \|_p ~=~ \inf_{{\bf z} \in \mathbb{C}^d, \| {\bf z} \|_0 \leq s} \left\| {\bf x} - {\bf z} \right\|_p
\end{equation}
for all $p \geq 1$.

\subsection{Bounded Orthonormal Systems}
\label{sec:BONS}

Let $\mathcal{D} \subset \mathbbm{R}^n$ be endowed with a probability measure $\mu$.  Further, let $\Psi = \{ \psi_0, \dots, \psi_{N} \}$ be an orthonormal set of real-valued functions on $\mathcal{D}$ so that
$$\int_{\mathcal{D}} \psi_i\left(~{\bf x} ~\right) \overline{\psi_j \left(~{\bf x}~ \right)} d \mu \left( {\bf x} \right) ~=~ \delta_{i,j}.$$
We will refer to any such $\Psi$ as an \textit{orthonormal system}.  More specifically, we will utilize a particular type of orthonormal system.
\begin{Def}
We call $\Psi = \{ \psi_0, \dots, \psi_N \}$ a bounded orthonormal system with constant $K \in \mathbbm{R}^+$ if
$$\left\| \psi_k \right\|_{\infty} := \sup_{{\bf x} ~\in~ \mathcal{D}} \left| \psi \left(~ {\bf x}~\right) \right| ~\leq~K~~\textrm{for all}~k \in [N].$$
\label{Def:BOS}
\end{Def}

For any orthonormal system, $\Psi$, on $\mathcal{D} \subset \mathbbm{R}^n$ with probability measure $\mu$, we may create an associated \textit{random sampling matrix}, $R \in \mathbbm{R}^{(m+1) \times (N+1)}$, as follows:  First, select $m+1$ points ${\bf x}_0, \dots, {\bf x}_m \in \mathcal{D}$ independently at random according to $\mu$.\footnote{So that $\mathbbm{P} \left[{\bf x}_j \in \mathcal{S} \right] = \mu \left( \mathcal{S} \right)$ for all measurable $\mathcal{S} \subseteq \mathcal{D}$ and $j \in [m]$.}  Then, form the matrix $R$ by setting $R_{i,j} := \psi_j \left(~ {\bf x}_i~\right)$ for each $i \in [m] \ $ and $j \in [N]$.  The following theorem concerning random sampling matrices created from bounded orthonormal systems in this fashion is proven in \cite{HolgerBook}.\footnote{See Theorem 12.31 in \cite{HolgerBook}.}

\begin{thm}
Let $R \in \mathbbm{R}^{(m+1) \times (N+1)}$ be a random sampling matrix created from a bounded orthonormal system with constant $K \geq 1$.  Let $\epsilon \in (0,1)$, $s \in [N] \setminus \{  0\}$, and set $\widetilde{R} = \frac{1}{\sqrt{m+1}}R$.  If $m \geq C K^2 \epsilon^{-2} s \ln^4(N)$, then with probability at least $1-N^{-\ln^3(N)}$
\begin{equation}
\sqrt{1 - \epsilon} ~\leq~ \sigma_{|\mathcal{S}|} \left(\widetilde{R}_{\mathcal{S}} \right)~\leq~\sigma_1\left( \widetilde{R}_{\mathcal{S}} \right) ~\leq~ \sqrt{1 + \epsilon}
\label{thm:randSampeq}
\end{equation}
will hold simultaneously for all nonempty subsets $\mathcal{S} \subset [N]$ having $|S| \leq s$ .  Here the constant $C > 0$ is fixed and universal.
\label{thm:randSamp}
\end{thm}

As pointed out in \cite{rauhut2012sparse}, the reweighed Legendre polynomials
\begin{equation}
Q_m(x) := \left(\frac{\pi}{2}\right)^{1/2} (1-x^2)^{1/4} L_m(x)
\label{equ:RewLeg}
\end{equation}
form a bounded orthonormal system with constant $K \leq \sqrt{3}$ with respect to the Chebyshev probability measure $d \mu(x) = \pi^{-1} (1 - x^2)^{-1/2} dx$ on $\mathcal{D} = [-1,1]$.  As a result, one easily obtains the following corollary of Theorem~\ref{thm:randSamp}.

\begin{cor}
Let $R \in \mathbbm{R}^{(m+1) \times (N+1)}$ be a random sampling matrix created from $\left\{ Q_0, \dots, Q_N \right\}$ in \eqref{equ:RewLeg} via $m+1$ points $x_0, \dots, x_m \in [-1,1]$ drawn independently according to the Chebyshev measure.  Let $\epsilon \in (0,1)$, $s \in [N] \setminus \{  0\}$, and set $\widetilde{R} = \frac{1}{\sqrt{m+1}}R$.  If $m \geq 3 C \epsilon^{-2} s \ln^4(N)$, then with probability at least $1-N^{-\ln^3(N)}$
$$\kappa \left(\widetilde{R}^*_{\mathcal{S}}\widetilde{R}_{\mathcal{S}} \right) ~\leq~ \frac{1 + \epsilon}{1 - \epsilon}$$
will hold simultaneously for all nonempty subsets $\mathcal{S} \subset [N]$ having $|S| \leq s$ .  The constant $C > 0$ is as in Theorem~\ref{thm:randSamp}.
\label{coro:randSamp}
\end{cor}

Corollary~\ref{coro:randSamp} guarantees that we can select a set of $m+1$ points from $[-1,1]$ \textit{once for each value of $N$} which will lead to a random sampling matrix for \eqref{equ:RewLeg}, $R \in \mathbbm{R}^{(m+1) \times (N+1)}$, all of whose associated submatrices $\widetilde{R}_{\mathcal{S}}$ are nearly isometric maps from $\mathbbm{R}^{|S|}$ into $\mathbbm{R}^m$.  Furthermore, $R$ can be formed in $\mathcal{O}(mN)$-time by using the standard recurrence relation for Legendre polynomials \cite{zhang1996computation} $$(m+1)L_{m+1}(x)=(2m+1)xL_m(x)-mL_{m-1}(x)$$
in order to generate each row.  This is a one-time computational cost for each choice of $m, N \in \mathbbm{N}$.  Alternatively, in a low memory setting, one may use fast methods based asymptotic expansions in order to quickly generate any desired submatrix of $\widetilde{R}$ from Corollary~\ref{coro:randSamp}, $\widetilde{R}_{\mathcal{S}} \in \mathbbm{R}^{(m+1) \times s}$, on the fly.  This can be accomplished in $\mathcal{O}(m \cdot s)$-time for any particular such submatrix of $\widetilde{R}$ as needed \cite{bogaert20121}.

\subsection{Iserles' Map from Fourier to Legendre Coefficients}
\label{sec:IserlesMap}

For a given analytic function $f$ and $r \in (0,1]$, define $f_r: \mathbb{C} \rightarrow \mathbb{C}$ by
\begin{equation}
f_r( x ) := \left( 1 - r^2 \mathbbm{e}^{2 \mathbbm{i} x} \right) f\left( \frac{1}{2}\left( r^{-1} \mathbbm{e}^{- \mathbbm{i} x} + r \mathbbm{e}^{\mathbbm{i} x} \right) \right).
\label{def:Fmap}
\end{equation}
Here, when $r<1$, $f$ is evaluated on a Bernstein ellipse in the complex plane; when $r=1$, $f$ is composed with $\cos(x)$ as per Chebyshev interpolation.  Note that $f_r$ will be both analytic (since $f$ is), as well as $2 \pi$-periodic on $\mathbb{R}$.  Hence, both the Legendre coefficients of $f: [-1, 1] \rightarrow \mathbb{R}$ and the Fourier series coefficients of $f_r:  [-\pi, \pi] \rightarrow \mathbb{C}$ will decay exponentially (see, e.g., \cite{titchmarsh1939theory,wang2012convergence}).  In any such setting the following map may be constructed from $\widehat{f}_r$ to $\tilde{f}$ (see \cite{iserles2011fast} for details).

Let $(a)_j$ be defined recursively for all $a \in \mathbb{R}$ and $j \in \mathbb{N}$ by
\begin{equation}
(a)_{j} := (a)_{j-1} \left( a+j-1 \right),
\label{Def:PochhammerSym}
\end{equation}
where $(a)_0 := 1$, and set
\begin{equation}
\tilde{g}_{i,j} := \frac{2^{2i} (i!)^2 (i+1)_j (\frac{1}{2})_j}{(2i)! j! (i + \frac{3}{2})_j} \cdot r^{i+2j},
\label{def:gmj}
\end{equation}
for all $i, j \in \mathbb{N}$.  Then, we have
\begin{equation}
\tilde{f}(i) = \sum^{\infty}_{j=0} \tilde{g}_{i,j} ~\widehat{f}_r \left( -i-2j \right)
\label{def:IMap}
\end{equation}
for all $i \in \mathbb{N}$.  Given the rapid decay of both $\tilde{g}_{i,j}$ and $\widehat{f}_r$ when $r \in (0,1)$, one may truncate the sum in order to approximate the first $N$ Legendre coefficients using
\begin{equation}
\tilde{f}(i) \approx \sum^{M}_{j=0} \tilde{g}_{i,j} ~\widehat{f}_r \left( -i-2j \right),
\label{equ:IserlesOriginal}
\end{equation}
for a modest $M \sim \log_{1/r}(N)$, after approximating the Fourier coefficients $\widehat{f}_r \left( 0 \right), \dots, \widehat{f}_r \left( -N-2M \right)$ using an FFT.  For $r = 1$ (or close to 1) a modest $M$ can still be chosen based solely on the decay of $\widehat{f}_r$.  The resulting numerical method requires $\mathcal{O} (N \log N)$ floating point operations in order to approximate $\tilde{f}(i)$ for all $i \in [N]$.

Herein we are primarily concerned with the setting where $f(x)$ is a polynomial of degree at most $N$ (recall $\eqref{Probf}$).  In this case it is easy to verify that $\widehat{f}_r(\omega) = 0$ for all $\omega < -N$ and $r \in (0,1]$.  Thus, the map $\eqref{def:IMap}$ reduces to the finite linear system
\begin{equation}
{\bf \tilde{f}} := \left( \begin{array}{l}  \tilde{f}(0) \\ ~~\vdots \\  \tilde{f}(N) \end{array} \right) = \tilde{G}_r \left( \begin{array}{l}  \widehat{f}_r \left( 0 \right) \\ ~~~~\vdots \\ \widehat{f}_r \left( -N \right) \end{array} \right)
\label{equ:IserlesNew}
\end{equation}
where $\tilde{G}_r \in \mathbb{R}^{(N+1) \times (N+1)}$ is the upper triangular matrix with entires
\begin{equation}
(\tilde{G}_r)_{i,j} := \left \{ \begin{array}{ll} \tilde{g}_{i,(j-i)/2} & {\rm if}~i \leq j, ~{\rm and}~ i \equiv j ~{\rm mod}~2 \\ 0 & {\rm else} \end{array} \right..
\end{equation}
We are now prepared to describe our method for rapidly and accurately computing sparse Legendre coefficient expansions.

\section{A Simple SFT-based approach for Reconstructing Sparse Legendre Expansions}
\label{sec:Alg}

We propose a two stage method for approximating functions, $f: [-1,1] \rightarrow \mathbb{R}$, with sparse/ compressible Legendre coefficient expansions as per \eqref{Probf}.  During the first stage, we use Iserles' map from \S \ref{sec:IserlesMap} in order to help us identify Legendre polynomials whose coefficients are large in magnitude in $f$.  We accomplish this by sampling $f$ according to its modified form, $f_r$ from \eqref{def:Fmap}, in order to take advantage of the fact that $\widehat{f}_r \left( -j \right) \approx r^{-j} \tilde{f}(j) / \sqrt{\pi j}$ holds for all $j \in [N]$ (see, e.g., Lemma~\ref{lem:EntryBound} together with Theorem~\ref{thm:InverseRowDecay} in \S \ref{sec:ErrorAnal}).  This fact guarantees that the Fourier coefficients of $f_r$ will be compressible whenever the Legendre coefficients of $f$, ${\bf \tilde{f}}$ from \eqref{equ:IserlesNew}, are sparse (see, e.g., Lemma~\ref{lem:FourierCompressible} in \S \ref{sec:ErrorAnal} for details).  Hence, we may utilize SFT techniques from \S \ref{sec:SFT} in order to rapidly identify the largest magnitude Fourier coefficients of $f_r$ which, in turn, immediately reveal the largest magnitude Legendre coefficients of $f$ via \eqref{equ:IserlesNew}.

We are happy:  using Iserles' map with an SFT is good enough to guarantee that one can rapidly identify large magnitude Fourier coefficients of $f_r$ which generally correspond to large magnitude Legendre coefficients of $f$!  Unfortunately, this SFT-based technique does not appear to allow us to quickly compute the Legendre coefficients with much accuracy.  SFTs can get us fast estimates that are ``in the right ball park" -- accurate enough to tell us that a coefficient is large -- but getting more than a few digits of accuracy this way appears elusive.  
Of course, one can always force a good SFT method to supply Iserles' method \eqref{equ:IserlesOriginal} with enough Fourier coefficients of $f_r$ to make it produce accurate estimates.  However, this appears to be far too slow an approach to be terribly interesting for the values of $r$ that produce decently conditioned maps $\tilde{G}_r$ (i.e., for $r \approx 1$).  

Thankfully, bounded orthonormal system results for reweighted Legendre polynomials (recall \S \ref{sec:BONS}) can easily solve the Legendre coefficient estimation problem for us {\it once we know which coefficients to compute}.  Given access to a 
well-conditioned random sampling matrix $R$ (recall Corollary~\ref{coro:randSamp}), along with small vector ${\bf y} \in \mathbbm{R}^{\mathcal{O}(s \log^4 N)}$ of additional reweighted samples from $f$ taken at the points $x_j \in [-1,1]$ used to build $R$,
\begin{equation}
y_j  = \left( \widetilde{R} {\bf \tilde{f}} \right)_j = \frac{\sqrt{\pi}(1-x_j^2)^{1/4}}{\sqrt{2(m+1)}} f(x_j),
\label{equ:ysamps}
\end{equation}
one can accurately estimate any given set of Legendre coefficients of $f$, $S \subset [N]$ with $|S| \leq s$, by quickly solving a small least-squares problem.  In particular, this means that we can simply $(i)$ identify a set of important Legendre coefficients of $f$ with an SFT, and then $(ii)$ use a random sampling matrix to accurately estimate the identified Legendre coefficients.  See Algorithm~\ref{alg:SFTLegendre} for pseudocode, and Lemma~\ref{lem:CoefEstLemma} in \S \ref{sec:ErrorAnal} for more details regarding coefficient estimation.

\begin{algorithm}
\renewcommand{\algorithmicrequire}{\textbf{Input:}}
\renewcommand{\algorithmicensure}{\textbf{Output:}}
\caption{Fast Sparse Legendre Coefficient Expansion Algorithm}
\label{alg:SFTLegendre}
\begin{algorithmic}[1]
    \REQUIRE $(i)$ Sparsity $s \in [N]$, $(ii)$ $r \in (0,1]$, $(iii)$ Pointer to a Random Sampling Matrix, $R \in \mathbbm{R}^{(m+1) \times (N+1)}$, as per Corollary~\ref{coro:randSamp}, $(iv)$ Renormalized samples from $f$, ${\bf y} \in \mathbbm{R}^{m+1}$, as per \eqref{equ:ysamps}, $(v)$ Pointer to $f$ from \eqref{Probf}, and $(vi)$ Pointer to a Sparse Fourier Transform code, ${\rm SFT}$
    \ENSURE A Sparse Approximation of the Legendre Coefficients of $f$ from \eqref{Probf}, ${\bf \tilde{f}'}_{\mathcal{S}}$
    \STATE Find the degrees of important Legendre polynomials present in $f$, $\mathcal{S} \subset [N]$ with $|\mathcal{S}| \leq s$, by running an SFT on $f_r$ from \eqref{def:Fmap} and recording the $\leq s$ most energetic frequencies it returns.
    \STATE Approximately Solve a Least Squares Problem:  ${\bf \tilde{f}'}_{\mathcal{S}} \approx {\bf z}_{\rm min} := \argmin_{{\bf z} \in \mathbbm{R}^{|\mathcal{S}|}} \left \| \widetilde{R}_{\mathcal{S}} {\bf z} - {\bf y} \right \|_2$.
\end{algorithmic}
\label{alg:FSL}
\end{algorithm}

The runtime complexity of Algorithm~\ref{alg:FSL} will generally be largely determined by the type of SFT chosen in line 1.  Both randomized and deterministic algorithms exist.  The randomized approaches are generally faster, but have a small (usually tunable) probability of failing to return a good answer.  The deterministic approaches are slower, but are guaranteed to approximate a given function as well as is possible with a sparse representation of the chosen size.  Recall Theorem~\ref{thm:MultDimRecov} in \S \ref{sec:SFT} for example results.

Considering the runtime complexity of line 2, we note that we may efficiently solve the least squares problem there using a Conjugate Gradient (CG) algorithm.  Suppose that the normalized random sampling matrix, $\widetilde{R} \in \mathbbm{R}^{(m+1) \times (N+1)}$, passed to Algorithm~\ref{alg:FSL} satisfies \eqref{thm:randSampeq} of Theorem~\ref{thm:randSamp} with $\epsilon = 3/5$.  In this case, a CG method will allow one to use just
$$C \log_{\frac{\sqrt{\kappa \left(\widetilde{R}^*_{\mathcal{S}}\widetilde{R}_{\mathcal{S}} \right)} + 1}{\sqrt{\kappa \left(\widetilde{R}^*_{\mathcal{S}}\widetilde{R}_{\mathcal{S}} \right)} - 1}} \left( \frac{\| {\bf y} \|_2}{\delta} \right) \leq C \log_{\frac{\sqrt{\frac{1 + \epsilon}{1 - \epsilon}} + 1}{\sqrt{\frac{1 + \epsilon}{1 - \epsilon}} - 1}} \left( \frac{\| {\bf y} \|_2}{\delta} \right) = C \log_{3} \left( \frac{\| {\bf y} \|_2}{\delta} \right)$$
CG iterations in order to get 
\begin{equation}
\left \| \widetilde{R}_{\mathcal{S}} \left({\bf z}_{\rm min} - {\bf \tilde{f}'}_{\mathcal{S}} \right) \right \|_2 \leq \delta.
\label{equ:ApproxCGsolve}
\end{equation}
for any desired $\delta \in \mathbbm{R}^+$ (see, e.g., Chapter 7 of \cite{bjorck1996numerical}).  Here, Corollary~\ref{coro:randSamp} (i.e., \eqref{thm:randSampeq}) has been used to bound $\kappa \left(\widetilde{R}^*_{\mathcal{S}}\widetilde{R}_{\mathcal{S}} \right)$ under the assumption that $m = C' s \ln^4(N)$ for appropriate fixed universal constants $C,C' \in \mathbbm{R}^+$.  Each CG iteration then takes $\mathcal{O} \left(s^2 \ln^4(N) \right)$-time.  Noting that $\| {\bf y} \|_2 \leq C'' \left \| {\bf \tilde{f}} \right \|_1$ will also hold, for a universal constant $C'' \in \mathbbm{R}^+$, whenever $\widetilde{R}$ satisfies \eqref{thm:randSampeq} (see, e.g., Exercise 6.6 in \cite{HolgerBook}), we have that we can compute an ${\bf \tilde{f}'}_{\mathcal{S}} \in \mathbbm{R}^{|\mathcal{S}|}$ satisfying \eqref{equ:ApproxCGsolve} in $\mathcal{O} \left(s^2 \ln^4(N) \cdot \ln \left( \frac{\left \| {\bf \tilde{f}} \right \|_1}{\delta} \right) \right)$-time.

\section{Error Analysis and Recovery Guarantees}
\label{sec:ErrorAnal}

In this section we analyze Algorithm~\ref{alg:SFTLegendre}.  The main results establish both that $(i)$ the largest magnitude Legendre coefficients present in $f$ can be rapidly identified via SFT methods (see Theorem~\ref{thm:SupportIdentification}), and that $(ii)$ once identified, the largest magnitude Legendre coefficients can be both rapidly and accurately approximated (see Lemma~\ref{lem:CoefEstLemma}).  By combining these results one can establish deterministic\footnote{Note that we are implicitly using randomized techniques to construct the random sampling matrices, $R \in \mathbbm{R}^{(m+1) \times (N+1)}$, used in line 2 of Algorithm~\ref{alg:FSL}.  However, the related probabilistic guarantees establish results \textit{for all sufficiently sparse signals} with high probability, and so can be viewed as establishing the existence of entirely deterministic methods.} sublinear-time recovery guarantees for many different classes of Legendre-compressible functions.  For example, one can easily prove sublinear-time recovery guarantees for exactly $s$-sparse Legendre polynomials of the form
\begin{equation}
f(x) = \sum_{n \in \mathcal{S} \subset [N]} \tilde{f}(n) L_n(x)
\label{equ:ExampPoly}
\end{equation}
where $|\mathcal{S}| = s$, $\min \mathcal{S} = \Omega(N/s)$, and all $s$ nonzero $\tilde{f}(n) \in \mathbbm{R}$ have (roughly) the same magnitude.  Doing so we may obtain, e.g., Theorem~\ref{thm:ExactSparseRecov}.

\begin{thm}
There exists a deterministic $\mathcal{O} \left( s^6 \log^5(N) \right)$-time algorithm that is guaranteed to exactly recover (up to machine precision) the Legendre coefficients of any function of type  \eqref{equ:ExampPoly}.
\label{thm:ExactSparseRecov}
\end{thm}

\noindent \textit{Proof:}  Apply Corollary~\ref{cor:ExactSparse} followed by Lemma~\ref{lem:CoefEstLemma}. \qed\\

Note that the runtime of the deterministic algorithm referred to by Theorem~\ref{thm:ExactSparseRecov} is indeed \textit{sub-linear in $N$ for sparsities $s \ll N$}.  However, it is also almost certainly suboptimal -- algorithmic modifications can probably be made that reduce the runtime complexity further without negatively impacting the recovery guarantee.  

It is also important to point out that the current assumptions concerning $f$ in \eqref{equ:ExampPoly} can be loosened considerably, without loosing deterministic recovery guarantees, by applying the subsequent results differently than done to get Theorem~\ref{thm:ExactSparseRecov}.  However, the theoretical results thus derived suffer both aesthetically and, in other ways, technically.   For this reason we will leave the proof of alternative guarantees via Theorem~\ref{thm:SupportIdentification} and Lemma~\ref{lem:CoefEstLemma} to the interested reader.  

We will now begin to prove our main theoretical results, starting with those concerning the rapid identification of the Legendre polynomials whose coefficients are largest in magnitude in $f$.

\subsection{Support Identification}

For the purposes of analyzing line 1 of Algorithm~\ref{alg:FSL} it is crucial to understand how sparse 
$${\bf \widehat{f}'_r} := \left( \begin{array}{l}  \widehat{f}_r \left( 0 \right) \\ ~~~~\vdots \\ \widehat{f}_r \left( -N \right) \end{array} \right)$$
will be given that ${\bf \tilde{f}}$ is sparse (recall \eqref{equ:IserlesNew}).  We will begin to move toward this goal by considering the matrix $\tilde{G}^{-1}_r$.  Once it is properly understood, we will then be able to consider the compressibility characteristics of ${\bf \widehat{f}'_r} = \tilde{G}^{-1}_r {\bf \tilde{f}}$ for sparse vectors ${\bf \tilde{f}}$.  The following lemma gives the entries of $\tilde{G}^{-1}_r$.

\begin{lem}
The even rows of the inverse matrix $\tilde{G}^{-1}_r$ from \eqref{equ:IserlesNew} are given by
\begin{equation}\label{EvenrowsofInverse}
\left( \tilde{G}^{-1}_r \right)_{2i,2j} = \left\{ \begin{array}{ll} 0 & j < i \\ \\
\frac{(-1)^j}{4^j} 
\sum^j_{k=i}  \frac{(-1)^k}{4^k} \frac{ r^{-2i} (2j+2k)!}{(j-k)! (j + k)! (2k)!} {2k \choose k+i} \frac{1+2i}{k+i+1}
& i \leq j \leq \left \lfloor \frac{N}{2} \right\rfloor\\ \\
 \end{array} \right.
\end{equation}
and $\left( \tilde{G}^{-1}_r \right)_{2i,2j+1} = 0$ for $i, j = 0, \dots, \left \lfloor \frac{N}{2} \right\rfloor$.  Similarly, the odd rows of the inverse matrix $\tilde{G}^{-1}_r$ from \eqref{equ:IserlesNew} are given by
\begin{equation}\label{OddrowsofInverse}
\left( \tilde{G}^{-1}_r \right)_{2i+1,2j+1} = \left\{ \begin{array}{ll} 0 & j < i \\ \\
\begin{array}{l}
\frac{(-1)^j}{2\cdot4^j} 
 \sum^j_{k=i}  \frac{(-1)^k}{2\cdot4^k} \frac{ r^{-2i-1} (2j+2k+2)!}{(j-k)! (j + k+1)! (2k+1)!} {2k+1 \choose k-i} \frac{2+2i}{k+i+2} 
 \end{array}&i \leq j <\left \lfloor \frac{N}{2} \right\rfloor\\ \\

 \end{array} \right.
\end{equation}
and $\left( \tilde{G}^{-1}_r \right)_{2i+1,2j} = 0$ for $i, j = 0, \dots, \left \lfloor \frac{N}{2} \right\rfloor$. 
\label{lem:InverseMatrix}
\end{lem}

\noindent \textit{Proof:}  See Appendix~\ref{app:ProofLemMatrixInv}. \qed\\

The following corollary of Lemma~\ref{lem:InverseMatrix} establishes simpler and more useful formulas for each entry of $\tilde{G}^{-1}_r$.

\begin{cor}
The even rows of the inverse matrix $\tilde{G}^{-1}_r$ from \eqref{EvenrowsofInverse} may also be written more compactly as
\begin{equation}\label{EvenrowsInv}
\left( \tilde{G}^{-1}_r \right)_{2i,2j} = \left\{ \begin{array}{ll} 0 & j < i \\ \\
 \left(-\frac{1}{4}\right)^{i+j}\frac{(2i+2j)!}{(i+j)!~(i+j)!} \frac{r^{-2i}{(2i+1)}}{i+j+1} \frac{\Gamma\left(\frac{3}{2}\right)}{(j-i)! ~ \Gamma\left(i-j+\frac{3}{2}\right)}
& i \leq j \leq \left \lfloor \frac{N}{2} \right\rfloor\\ \\
 \end{array} \right.
\end{equation}
and $\left( \tilde{G}^{-1}_r \right)_{2i,2j+1} = 0$ for $i, j = 0, \dots, \left \lfloor \frac{N}{2} \right\rfloor$.  Similarly, the odd rows of the inverse matrix $\tilde{G}^{-1}_r$ from \eqref{OddrowsofInverse} may be written as 
\begin{equation}\label{OddrowsInv}
\left( \tilde{G}^{-1}_r \right)_{2i+1,2j+1} = \left\{ \begin{array}{ll} 0 & j < i \\ \\
 \frac{(-1)^{i+j}}{4^{i+j+1}}\frac{(2(i+j+1))!}{(i+j+1)!~(i+j+1)!} \frac{r^{-2i-1}{(2i+2)}}{i+j+2}\frac{\Gamma\left(\frac{3}{2}\right)}{(j-i)! ~ \Gamma\left(i-j+\frac{3}{2}\right)}
& i \leq j \leq \left \lfloor \frac{N}{2} \right\rfloor\\ \\
 \end{array} \right.
\end{equation}
and $\left( \tilde{G}^{-1}_r \right)_{2i+1,2j} = 0$ for $i, j = 0, \dots, \left \lfloor \frac{N}{2} \right\rfloor$. 
\label{cor:CompactInverseMatrix}
\end{cor}

\noindent \textit{Proof:}  See Appendix~\ref{app:ProofLemCompactMatrixInv}. \qed\\

With Corollary~\ref{cor:CompactInverseMatrix} in hand, one may now see that each row and column of $\tilde{G}_{r}^{-1}$ is weakly dominated by it's diagonal entry. See Theorem~\ref{thm:InverseRowDecay} below for an exact statement.

\begin{thm}
For nonzero entries of the inverse matrix $\tilde{G}_{r}^{-1}$, the decay rate of each row is given by
\begin{equation}\label{rowdecay}
\begin{array}{rl}
\left| (\tilde{G}_{r}^{-1})_{n, n+2x}\right| <  \frac{\sqrt[3]{\textrm{e}}}{2\sqrt{\pi}} \frac{1}{x\sqrt{x-2}}\left| (\tilde{G}_{r}^{-1})_{n, n}\right| & \textrm{ for } x>2, \\
\left| (\tilde{G}_{r}^{-1})_{n, n+2} \right| < \frac{1}{2}\left| (\tilde{G}_{r}^{-1})_{n, n} \right|, \textrm{ and } \left| (\tilde{G}_{r}^{-1})_{n, n+4} \right| < \frac{1}{8}\left| (\tilde{G}_{r}^{-1})_{n, n} \right| & \textrm{ for } x=1, x=2. 
\end{array}
\end{equation}

The decay rate of each column is given by
\begin{equation}\label{columndecay}
\begin{array}{rl}
\left| (\tilde{G}_{r}^{-1})_{n-2x, n}\right| \le  \frac{\sqrt[3]{\textrm{e}}}{\sqrt{2\pi}}\frac{r^{2x}}{x\sqrt{x-2}} \left| (\tilde{G}_{r}^{-1})_{n, n}\right| & \textrm{ for } 2<x\le \left \lfloor \frac{n}{2} \right\rfloor, ~n \geq 6\\
\left| (\tilde{G}_{r}^{-1})_{n-2, n} \right| < \frac{r^2}{2}\left| (\tilde{G}_{r}^{-1})_{n, n} \right|, \textrm{and} \left| (\tilde{G}_{r}^{-1})_{n-4, n} \right| < \frac{r^4}{4} \left| (\tilde{G}_{r}^{-1})_{n, n} \right| & \textrm{ for } x=1, x=2, n \geq 2x. 
\end{array}
\end{equation}

The diagonal entries of $\tilde{G}_{r}^{-1}$ satisfy 
\begin{align}\label{diagonalentrydecay}
\frac{r^{-n}}{\sqrt{\pi n}} \left( 1-\frac{1}{8n}\right) \le \left| (\tilde{G}_{r}^{-1})_{n,n} \right| \le \frac{r^{-n}}{\sqrt{\pi n}} \, \textrm{ for } n > 0,~{\rm and}~\left| (\tilde{G}_{r}^{-1})_{0,0} \right| = 1.
\end{align}
\label{thm:InverseRowDecay}
\end{thm}

\noindent \textit{Proof:}  First, let's simplify the term $\Gamma\left(\frac{3}{2}\right)/\left((j-i)! \Gamma\left(i-j+\frac{3}{2}\right)\right)$. 
Since $\Gamma(3/2)=\sqrt{\pi}/2$, and by the reflection formula for the Gamma function (see, e.g., \cite{zhang1996computation}), $
\Gamma(-z)=-\pi/\left(z\Gamma(z)\sin(\pi z)\right)$ for $z\not\in \mathbb{Z}$, 
we have
\begin{align}
\frac{\Gamma\left(\frac{3}{2}\right)}{(j-i)! \Gamma\left(i-j+\frac{3}{2}\right)}=& \frac{(-1)^{j-i+1}}{2\sqrt{\pi}}\frac{(j-i-\frac{3}{2})\Gamma(j-i-\frac{3}{2})}{(j-i)!}. \nonumber
\end{align}
According to the half integer argument for the Gamma function, $\Gamma(n/2)=\left((n-2)!!\sqrt{\pi}\right)/2^{(n-1)/2}$ for $n\in \mathbb{Z}^{+}$, where $n!!=(2k)!/\left(2^{k}k!\right)$ when $n=2k-1$ for $k\in \mathbb{Z}^{+}$.  Thus,
\begin{align}
\frac{\Gamma\left(\frac{3}{2}\right)}{(j-i)! \Gamma\left(i-j+\frac{3}{2}\right)}=& \frac{(-1)^{j-i+1}}{2}\frac{(j-i-\frac{3}{2})}{(j-i)!}\frac{(2j-2i-5)!!}{2^{j-i-2}}\nonumber \\
=& \frac{(-1)^{j-i+1}}{2\cdot 4^{j-i-2}}\frac{(j-i-\frac{3}{2})}{(j-i)(j-i-1)}\frac{\left(2(j-i-2)\right)!}{(j-i-2)!(j-i-2)!}.\nonumber
\end{align}
By Stirling's approximation, $n!=\sqrt{2\pi n}n^{n}\textrm{e}^{-n}\textrm{e}^{R_{n}}$, where $0\le R_n \le 1/12n$ (see, e.g., \cite{robbins1955remark, HolgerBook}), 
\begin{align}\label{simpleform}
\frac{\Gamma\left(\frac{3}{2}\right)}{(j-i)! \Gamma\left(i-j+\frac{3}{2}\right)}=&\frac{(-1)^{j-i+1}}{2\sqrt{\pi}}\frac{(j-i-\frac{3}{2})}{(j-i)(j-i-1)\sqrt{j-i-2}}\textrm{e}^{R_{2(j-i-2)}-2R_{j-i-2}}.
\end{align}

Now by Corollary \ref{cor:CompactInverseMatrix} and \eqref{simpleform}, if $n=2i$ for $i,x \in \mathbb{Z}$, $i\ge 0$, and $x > 2$, then 
\begin{equation}\label{rowdecay_sim}
\begin{array}{ll}
\left| (\tilde{G}_{r}^{-1})_{n, n+2x} \right| 
&= \left| (\tilde{G}_{r}^{-1})_{n, n}\right|\cdot\left| \frac{1}{4^x} \frac{(n!)^2}{(2n)!}\frac{(2n+2x)!}{(n+x)!(n+x)!}\frac{n+1}{n+x+1}\frac{1}{2\sqrt{\pi}}\frac{x-\frac{3}{2}}{x(x-1)\sqrt{x-2}}\textrm{e}^{R_{2(x-2)}-2R_{x-2}}\right| 
\end{array}
\end{equation}
Using Stirling's approximation again for $i\ge 1$ we have 
\begin{eqnarray}
\left| (\tilde{G}_{r}^{-1})_{n, n+2x} \right| =& \left| (\tilde{G}_{r}^{-1})_{n, n}\right| \left( \frac{\sqrt{n}}{\sqrt{n+x}}\frac{n+1}{n+x+1}\frac{(x-\frac{3}{2})}{x-1}\frac{1}{x\sqrt{x-2}}\frac{1}{2\sqrt{\pi}} \textrm{e}^{R^{\textrm{a}}} \right), \label{equ:RelRowBound}
\end{eqnarray}
where $R^{\textrm{a}}=2R_{n}-R_{2n}+R_{2n+2x}-2R_{n+x}+R_{2(x-2)}-2R_{x-2} \leq 1/3$ for all $n,x \in \mathbbm{Z}^+$ with $x \geq 3$. 
According to \eqref{rowdecay_sim}, when $i=0$, it is easy to show that 
\begin{equation}
\left|(\tilde{G}_{r}^{-1})_{0, 2x}\right| \leq \left|(\tilde{G}_{r}^{-1})_{0, 0}\right|\cdot \frac{\sqrt[6]{\textrm{e}}}{2\pi}\frac{1}{(x+1)}\frac{1}{x^{3/2} \sqrt{x-2}},\, \textrm{ for }x>2.
\end{equation}
Similarly, if $n=2i+1$ for $i\ge 0$ it is not difficult to verify that \eqref{equ:RelRowBound} still holds.
Therefore, one can see that each row satisfies
\begin{eqnarray}
\left| (\tilde{G}_{r}^{-1})_{n, n+2x}\right| <  \frac{\sqrt[3]{\textrm{e}}}{2\sqrt{\pi}} \frac{1}{x\sqrt{x-2}}\left| (\tilde{G}_{r}^{-1})_{n, n}\right|. \nonumber
\end{eqnarray}
The remainder of the proof of \eqref{rowdecay} is now easily established using Corollary \ref{cor:CompactInverseMatrix}. 

Similarly, we can bound the rate of decay of each column off of the diagonal.  By Corollary \ref{cor:CompactInverseMatrix}, \eqref{simpleform}, and Stirling's approximation, if $n=2i$ for $i\ge 3$, $2<x\le \left \lfloor \frac{n}{2} \right\rfloor$, then
\begin{equation*}
\begin{array}{ll}
\left| (\tilde{G}_{r}^{-1})_{n-2x, n} \right| 
&= \left| (\tilde{G}_{r}^{-1})_{n, n}\right| \cdot \frac{r^{2x}}{2\sqrt{\pi}} \left( \sqrt{\frac{n}{n-x}} \right) \frac{n-2x+1}{n-x+1} \left( \frac{x-\frac{3}{2}}{x-1} \right) \frac{\textrm{e}^{\tilde{R}^{\textrm{a}}}}{x\sqrt{x-2}}
\end{array}
\end{equation*}
where $\tilde{R}^{\textrm{a}}=2R_{n}-R_{2n}+R_{2(n-x)}-2R_{n-x}+R_{2(x-2)}-2R_{x-2} \leq 1/3$ for all $n,x \in \mathbbm{Z}^+$ with $x \geq 3$.  Since $\sqrt{n/(n-x)}\le \sqrt{2}$ for $2<x\leq i$, 
\begin{equation}
\begin{array}{ll}
\left| (\tilde{G}_{r}^{-1})_{n-2x, n} \right| 
&< \left| (\tilde{G}_{r}^{-1})_{n, n}\right| \cdot  \frac{r^{2x}}{\sqrt{2\pi}}\frac{\sqrt[3]{\textrm{e}}}{x\sqrt{x-2}}. 
\end{array}
\label{equ:OddDecayColumn}
\end{equation}
For $n=2i+1$ an analogous calculation reveals that \eqref{equ:OddDecayColumn} still holds.  The remainder of the proof of \eqref{columndecay} is now easily established using Corollary \ref{cor:CompactInverseMatrix}.  Lastly, the proof of \eqref{diagonalentrydecay} is easily established using Lemma \ref{lem:InverseMatrix} together with Theorem 2.6 in \cite{stanica2001good}.  \qed\\

We are now in the position to begin studying the compressibility of ${\bf \widehat{f}'_r} = \tilde{G}^{-1}_r {\bf \tilde{f}}$ in terms of the compressibility of ${\bf \tilde{f}}$.  Let $\sigma: [N] \rightarrow [N]$ be a permutation of $[N]$ such that
$$\left|  \tilde{f}_{\sigma(j)} \right| \geq \left|  \tilde{f}_{\sigma(j+1)} \right|$$
holds for all $j \in [N-1]$.  Let $\rho:\mathbbm{R} \rightarrow \mathbbm{R}^+ \cup \{ \infty \}$ be a modified ramp function with $\rho(x) = x$ for all $x \in \mathbbm{R}^+$, and $\rho(x) = \infty$ for all $x < 0$.  
Finally, define the right-distance form $j \in [N]$ to the set $\left\{ \sigma(n) ~\Big|~ n \in [s-1],~\sigma(n) \equiv j ~{\rm mod}~2 \right \}$ to be 
$$d_s(j) := \min \left( \left \{ \frac{\rho \left( \sigma(n) - j \right)}{2} ~\Big|~ n \in [s-1],~\sigma(n) \equiv j ~{\rm mod}~ 2 \right\} \cup \{ \infty \} \right).$$
We now have sufficient notation to consider the sizes of specific entries of ${\bf \widehat{f}'_r}$.

\begin{lem}
Let $s \in [N] \setminus \{0 \}$.  We have that
$$\left| \left( \widehat{f}'_r \right)_j -  (\tilde{G}^{-1}_r)_{j,j} \tilde{f}_j \right|  < \frac{31}{20} \left| (\tilde{G}^{-1}_r)_{j,j} \right| \left| \tilde{f}_{\sigma(s)} \right| +  \frac{\sqrt[3]{\textrm{e}}}{2\sqrt{\pi}} \left( \frac{ \| {\bf \tilde{f}} \|_1}{d_s(j)\sqrt{d_s(j)-2}} \right)  \left| (\tilde{G}^{-1}_r)_{j,j} \right|$$
holds for all $j \in [N]$ with $d_s(j) > 2$.  
\label{lem:EntryBound}
\end{lem}

\noindent \textit{Proof:}  Using that ${\bf \widehat{f}'_r} = \tilde{G}^{-1}_r {\bf \tilde{f}}$ with $\tilde{G}^{-1}_r$ upper triangular as per Lemma~\ref{lem:InverseMatrix}, we have that
$$\left( \widehat{f}'_r \right)_j =  (\tilde{G}^{-1}_r)_{j,j} \tilde{f}_j + \sum^{\left \lfloor \frac{N-j}{2} \right \rfloor}_{l = 1} (\tilde{G}^{-1}_r)_{j,j+2l} \tilde{f}_{j+2l}.$$
Rearranging this expression we can see that
\begin{align*}
\left| \left( \widehat{f}'_r \right)_j -  (\tilde{G}^{-1}_r)_{j,j} \tilde{f}_j \right|  &\leq \sum^{d_s(j)-1}_{l = 1} \left| (\tilde{G}^{-1}_r)_{j,j+2l} \right| \left|  \tilde{f}_{j+2l} \right| + \sum^{\left \lfloor \frac{N-j}{2} \right \rfloor}_{l = d_s(j)} \left| (\tilde{G}^{-1}_r)_{j,j+2l} \right| \left|  \tilde{f}_{j+2l} \right| \\
& <  \left|  \tilde{f}_{\sigma(s)} \right|  \left( \sum^{d_s(j)-1}_{l = 1} \left| (\tilde{G}^{-1}_r)_{j,j+2l} \right| \right) +  \frac{\sqrt[3]{\textrm{e}}}{2\sqrt{\pi}} \left( \frac{ \| {\bf \tilde{f}} \|_1}{d_s(j)\sqrt{d_s(j)-2}} \right)  \left| (\tilde{G}^{-1}_r)_{j,j} \right|,
\end{align*}
where the second inequality follows from H\"older's Inequality, Theorem~\ref{thm:InverseRowDecay}, and the definition of $d_s(j)$.  Appealing to Theorem~\ref{thm:InverseRowDecay} again we can also see that
\begin{equation*}
\sum^{d_s(j)-1}_{l = 1} \left| (\tilde{G}^{-1}_r)_{j,j+2l} \right| < \left| (\tilde{G}^{-1}_r)_{j,j} \right| \left( \frac{1}{2} + \frac{1}{8} + \frac{\sqrt[3]{\textrm{e}}}{6\sqrt{\pi}} + \frac{\sqrt[3]{\textrm{e}}}{2\sqrt{\pi}} \int^{\infty}_3 \frac{dx}{x \sqrt{x-2}}\right).
\end{equation*}
The remainder of the proof now follows.    \qed\\

We can now begin to understand the compressibility of ${\bf \widehat{f}'_r}$.  In particular, we have the following result concerning the best $s \cdot k$ approximation to ${\bf \widehat{f}'_1}$ for any $k \in \mathbbm{Z}^+$ with $k > 5$.

\begin{lem}
Let $s,k \in [N]$ with $k > 5$, $s > 8$, and $n_{\rm min} := \min \left\{ \sigma(n) ~\Big|~ n \in [s-1] \right\} > k$.  Then, 
$$\left \| {\bf \widehat{f}'_1} - \left({\bf \widehat{f}'_1} \right)^{\rm~ opt}_{s \cdot k} \right\|_1 < 7 \sqrt{N} \left| \tilde{f}_{\sigma(s)} \right| + \frac{s \| {\bf \tilde{f}} \|_1}{\sqrt{n_{\rm min}-4} \sqrt{k-5}}.$$
\label{lem:FourierCompressible}
\end{lem}

\noindent \textit{Proof:}  We bound the sum of the magnitudes of all entries in ${\bf \widehat{f}'_1}$ whose indices have right-distance $\geq k$ from $\left\{ \sigma(n) ~\Big|~ n \in [s-1] \right \}$ using two cases:  

{\bf Case I.}  We bound the magnitudes of the small entries with $j \leq n_{\rm min} - k$.  Using Lemma~\ref{lem:EntryBound} we can see that
\begin{equation*}
\sum^{n_{\rm min} - k}_{j=0} \left| \left( \widehat{f}'_1 \right)_j \right| <  \sum^{n_{\rm min} - k}_{j=0} \left| (\tilde{G}^{-1}_r)_{j,j} \right| \left( \frac{51}{20} \left| \tilde{f}_{\sigma(s)} \right| + \frac{\sqrt[3]{\textrm{e}}}{\sqrt{\pi / 2}} \left( \frac{ \| {\bf \tilde{f}} \|_1}{(n_{\rm min}-j)\sqrt{n_{\rm min}-j-4}} \right)  \right) .
\end{equation*}
Theorem~\ref{thm:InverseRowDecay} now implies that 
\begin{align*}
\sum^{n_{\rm min} - k}_{j=0} \left| \left( \widehat{f}'_1 \right)_j \right| <~ & \frac{51 \sqrt{2 \cdot n_{\rm min} }}{10 \sqrt{\pi}} \left| \tilde{f}_{\sigma(s)} \right| +\frac{2 \| {\bf \tilde{f}} \|_1}{n_{\rm min} - 1}+ \frac{\sqrt[3]{\textrm{e}} \sqrt{2} \| {\bf \tilde{f}} \|_1}{\pi} \cdot \sum^{n_{\rm min} - k}_{j=2}  \frac{1}{\sqrt{j}(n_{\rm min}-j-4)^{\frac{3}{2}}} \\
<~ & \frac{51 \sqrt{2 \cdot n_{\rm min} }}{10 \sqrt{\pi}} \left| \tilde{f}_{\sigma(s)} \right| +\frac{2 \| {\bf \tilde{f}} \|_1}{n_{\rm min} - 1}+ \frac{\sqrt[3]{\textrm{e}} \sqrt{2} \| {\bf \tilde{f}} \|_1}{\pi} \cdot \int_1^{n_{\rm min} - k+1} \frac{dx}{\sqrt{x} \left( n_{\rm min}-4-x \right)^{3/2}} \\
<~ & 4.1 \sqrt{n_{\rm min} } \left| \tilde{f}_{\sigma(s)} \right| +\frac{4 \| {\bf \tilde{f}} \|_1}{\sqrt{n_{\rm min} - 4} \sqrt{k-5}} .
\end{align*}

{\bf Case II.}  Here we bound the sum of the magnitudes of all entries $\left| \left( \widehat{f}'_1 \right)_j \right|$ with $j \in A := \{ x ~|~ x > n_{\rm min}$ and $d_s(x) \geq k \}$.  Note that there can be at most $(s-1) k$ entries in $[N] \setminus (A \cup [n_{\rm min}])$.  Using Lemma~\ref{lem:EntryBound}, Theorem~\ref{thm:InverseRowDecay}, and definition of $A$ we can see that
\begin{align*}
\sum_{j \in A} \left| \left( \widehat{f}'_1 \right)_j \right| &~<~   \frac{51}{20 \sqrt{\pi}} \left| \tilde{f}_{\sigma(s)} \right| \left( \sum^{N}_{j=n_{\rm min} + 1} \frac{1}{\sqrt{j}} \right)+ \frac{\sqrt[3]{\textrm{e}} (s-1)}{2\pi \sqrt{n_{\rm min}}} \left( \sum^{\infty}_{l=k}  \frac{ \| {\bf \tilde{f}} \|_1}{l\sqrt{l-2}} \right) \\
&~<~ \frac{51 \sqrt{N}}{10 \sqrt{\pi}} \left| \tilde{f}_{\sigma(s)} \right| + \frac{\sqrt[3]{\textrm{e}} (s-1)}{\pi \sqrt{n_{\rm min}}} \frac{ \| {\bf \tilde{f}} \|_1}{\sqrt{k-3}}.
\end{align*}
Combining the bounds from Cases I and II now finishes the proof.  \qed\\

Note that the vector ${\bf \widehat{f}'_r}$ contains only the (potentially nonzero) negative Fourier series coefficients of $f_1(x) =  \left( 1 - \mathbbm{e}^{2 \mathbbm{i} x} \right) f\left( \cos(x) \right)$.  Let 
$${\bf \widehat{f}_1} := \left( \begin{array}{l} \widehat{f}_1 \left( N + 2 \right) \\ ~~~~\vdots \\ \widehat{f}_1 \left( 0 \right) \\ ~~~~\vdots \\ \widehat{f}_1 \left( -N \right) \end{array} \right) \in \mathbbm{R}^{2N+3}$$
consist of all the potentially nonzero Fourier series coefficients of $f_1(x)$.  Noting that $f\left( \cos(x) \right)$ is an even real-valued function, one can see that $\widehat{f}_1 \left( \omega \right) = - \widehat{f}_1 \left( 2-\omega \right)$ holds for all $\omega \in \mathbbm{Z}^+$ (with $\widehat{f}_1 \left( 1 \right) = 0$).  As a result, Lemmas~\ref{lem:EntryBound} and~\ref{lem:FourierCompressible}  trivially extend to ${\bf \widehat{f}_1}$.  Using this fact in combination with results from \cite{I2013Improved} finally allows us to prove the main result of this section.

\begin{thm}
Let $s,k \in [N]$ with $k > 5$, $s > 8$, $\mathcal{S} := \left\{ \sigma(n) ~\Big|~ n \in [s-1] \right \}$, and $n_{\rm min} := \min \mathcal{S} > k$.  Given $j \in \mathcal{S}$, define the new right-distance from $j$ to $\mathcal{S}$ to be 
$$d'_s(j) := \min \left( \left \{ \frac{\rho \left( \sigma(n) - j \right)}{2} ~\Big|~ n \in [s-1],~\sigma(n) \equiv j ~{\rm mod}~ 2 \right\} \setminus \{ 0 \} \cup \{ \infty \} \right) > 0 = d_s(j).$$
Then, the deterministic variant of the algorithm referred to by Theorem~\ref{thm:MultDimRecov} will recover all $j \in \mathcal{S}$ satisfying both
\begin{equation}
\left|  \tilde{f}_j  \right| \geq \left( \frac{224 \sqrt{\pi} N}{7 s \cdot k} + \frac{31}{20} \right) \left|  \tilde{f}_{\sigma(s)} \right| + \frac{17 \sqrt{N} \| {\bf \tilde{f}} \|_1}{2 k \sqrt{(n_{\rm min}-4)(k-5)}},
\label{equ:BigLegendreCoef}
\end{equation}
and $d'_s(j) \geq k$.  Its operation count will be
\begin{equation*}
\mathcal{O} \left( s^2 \cdot k^2 \cdot \log^4 N \right).
\end{equation*}
\label{thm:SupportIdentification}
\end{thm}

\noindent \textit{Proof:}  We consider the behavior of the deterministic variant of the algorithm referred to by Theorem~\ref{thm:MultDimRecov} when executed on $f_1(x) = \left( 1 - \mathbbm{e}^{2 \mathbbm{i} x} \right) f\left( \cos(x) \right)$ with sparsity parameter $2s$ and $\epsilon = 1/k$.  In the course of forming its output trigonometric polynomial, $y_s$, this algorithm is guaranteed to identify every $j \in [-N, N+2] \cap \mathbbm{Z}$ with the property that 
$$\left| \widehat{f}_1  \left(  j \right)  \right|> \frac{2 \left \| {\bf \widehat{f}_1} - \left({\bf \widehat{f}_1} \right)^{\rm~ opt}_{2s \cdot k} \right\|_1}{s \cdot k}$$
by Lemma 6 in \cite{I2013Improved}.

Suppose that $\tilde{f}_j$ satisfies \eqref{equ:BigLegendreCoef} and also has $d'_s(j) \geq k$.  A trivial variant of Lemma~\ref{lem:EntryBound} then implies that
$$\left| \widehat{f}_1  \left( -j \right)  \right| > \left| (\tilde{G}^{-1}_r)_{j,j} \right|  \left| \tilde{f}_j \right|  - \frac{31}{20} \left| (\tilde{G}^{-1}_r)_{j,j} \right| \left| \tilde{f}_{\sigma(s)} \right| -  \frac{\sqrt[3]{\textrm{e}}}{2\sqrt{\pi}} \left( \frac{ \| {\bf \tilde{f}} \|_1}{d'_s(j)\sqrt{d'_s(j)-2}} \right)  \left| (\tilde{G}^{-1}_r)_{j,j} \right|.$$
Using \eqref{equ:BigLegendreCoef} and that $d'_s(j) \geq k$ one can now see that 
$$\left|  \widehat{f}_1  \left( -j \right)  \right| > \left| (\tilde{G}^{-1}_r)_{j,j} \right| \left( \frac{28 \cdot 8 \sqrt{\pi} N}{7 s \cdot k}  \left|  \tilde{f}_{\sigma(s)} \right| + \frac{32 \sqrt{\pi N} \| {\bf \tilde{f}} \|_1}{7 k \sqrt{(n_{\rm min}-4)(k-5)}}  \right).$$
Theorem~\ref{thm:InverseRowDecay} followed by Lemma~\ref{lem:FourierCompressible} finally implies that 
\begin{align*}
\left|   \widehat{f}_1  \left( -j \right)   \right| &> \frac{28 \cdot \sqrt{N}}{s \cdot k}  \left|  \tilde{f}_{\sigma(s)} \right| + \frac{4 \| {\bf \tilde{f}} \|_1}{k \sqrt{(n_{\rm min}-4)(k-5)}} \\
&= \frac{2}{s \cdot k} \left( 14 \sqrt{N} \left| \tilde{f}_{\sigma(s)} \right| + \frac{2s \| {\bf \tilde{f}} \|_1}{\sqrt{n_{\rm min}-4} \sqrt{k-5}} \right) \\
&> \frac{2 \left \| {\bf \widehat{f}_1} - \left({\bf \widehat{f}_1} \right)^{\rm~ opt}_{2s \cdot k} \right\|_1}{s \cdot k}.
\end{align*}
This guarantees that $j$ will indeed be identified as claimed.  \qed\\

The final corollary of this section applies Theorem~\ref{thm:SupportIdentification} to the situation of exactly $s$-sparse vectors ${\bf \tilde{f}}$.  We have the following result:

\begin{cor}
Let $s \in [N]$ with $s > 8$, $\mathcal{S} := \left\{ \sigma(n) ~\Big|~ n \in [s-1] \right \}$, and $n_{\rm min} := \min \mathcal{S} >  \frac{N}{17s/2-5} + 4 > 17 s / 2$.  Furthermore, suppose that ${\bf \tilde{f}}$ is exactly $s$-sparse with $\left|  \tilde{f}_{\sigma(0)} \right| = \left|  \tilde{f}_{\sigma(s-1)} \right| > \tilde{f}_{\sigma(s)} = 0$.  Then, there exists a deterministic algorithm which will recover a set of cardinality $\mathcal{O} \left( s^3 \right)$ that is guaranteed to contain $\mathcal{S}$ as a subset.  Its operation count will be $\mathcal{O} \left( s^4 \cdot \log^4 N \right)$.
\label{cor:ExactSparse}
\end{cor}

\noindent \textit{Proof:}  Apply Theorem~\ref{thm:SupportIdentification} with $k = 17 s / 2$.  The deterministic variant of the algorithm referred to by Theorem~\ref{thm:MultDimRecov} will recover all $j \in \mathcal{S}$ satisfying both $d'_s(j) \geq 17 s / 2$, and
\begin{equation*}
\left|  \tilde{f}_j  \right| \geq \frac{\sqrt{N} \| {\bf \tilde{f}} \|_1}{s \sqrt{(n_{\rm min}-4)(17 s / 2-5)}}, ~{\rm where}~ \frac{\sqrt{N} \| {\bf \tilde{f}} \|_1}{s \sqrt{(n_{\rm min}-4)(17 s / 2-5)}} < \frac{\| {\bf \tilde{f}} \|_1}{s} = \left|  \tilde{f}_{\sigma(s-1)} \right|,
\end{equation*}
as a subset of a set $\mathcal{A}$ of cardinality $\mathcal{O}(s^2)$.  Returning all $j \in [N]$ whose right-distance to $\mathcal{A}$ is $< 17 s / 2$ provides a set of cardinality $\mathcal{O} \left( s^3 \right)$ that is guaranteed to contain $\mathcal{S}$ as a subset.  The result follows. \qed\\

We conclude this section by mentioning several obvious facts regarding Theorem~\ref{thm:SupportIdentification} and Corollary~\ref{cor:ExactSparse}.  First, we have no doubt that they can be improved in general.  The easiest way to do this is to utilize randomized SFT methods in place of the deterministic algorithm from Theorem~\ref{thm:MultDimRecov} in order to reduce the computational complexities involved.    It is also clear that the assumptions in Corollary~\ref{cor:ExactSparse} can be relaxed rather easily at the price of increased computational costs.  Less trivial improvements would probably revolve around developing better variants of $f_1(x)$ whose Fourier coefficients are less contaminated by Legendre coefficients that are ``far away" from the associated ones (i.e., $\tilde{G}^{-1}_r$ should be ``more diagonally dominant").  Alternatively, one might design efficient filtering schemes which gradually reveal the lower (generally more energetic) frequencies of $f_1(x)$ so that they do not overwhelm larger (generally less energetic) frequencies as we hunt for them.  However, we will leave such considerations for future work.  We are now ready to consider how accurately we can estimate the Legendre coefficients for the set of  supporting polynomials, $\mathcal{S} \subset [N]$, we identify in this section.

\subsection{Coefficient Estimation}

Estimating the Legendre coefficients for the polynomials that have been identified as present in $f$ is comparatively easy given all the previous work on bounded orthonormal systems.  We have the following result:

\begin{lem}
Suppose that the normalized random sampling matrix, $\widetilde{R} \in \mathbbm{R}^{(m+1) \times (N+1)}$, passed to Algorithm~\ref{alg:FSL} satisfies \eqref{thm:randSampeq} of Theorem~\ref{thm:randSamp} with $\epsilon = 3/5$.  Let $\delta \in \mathbbm{R}^+$, and $\mathcal{S} \subset [N]$ have cardinality $|\mathcal{S}| = s$.  Then line 2 of of Algorithm~\ref{alg:FSL} can produce an $s$-sparse vector, ${\bf z} \in \mathbbm{R}^{N+1}$, satisfying 
$$\left \| {\bf \tilde{f}} - {\bf z}  \right \|_2 ~\leq~ 5 \left \| {\bf \tilde{f}}_{\mathcal{S}^{\rm c}} \right \|_2 + \frac{ 4 \left \| {\bf \tilde{f}}_{\mathcal{S}^{\rm c}} \right \|_1}{\sqrt{s}} + \sqrt{\frac{5}{2}} \delta$$
in $\mathcal{O} \left(s^2 \ln^4(N) \cdot \ln \left( \frac{\left \| {\bf \tilde{f}} \right \|_1}{\delta} \right) \right)$-time.
\label{lem:CoefEstLemma}
\end{lem}

\noindent \textit{Proof:}  Let ${\bf \tilde{f}''}_{\mathcal{S}} \in \mathbbm{R}^{|\mathcal{S}|}$ be such that $ \widetilde{R}_{\mathcal{S}} {\bf \tilde{f}''}_{\mathcal{S}} = \widetilde{R} {\bf \tilde{f}}_{\mathcal{S}}$ (i.e., let ${\bf \tilde{f}''}_{\mathcal{S}}$ be ${\bf \tilde{f}}_{\mathcal{S}}$ with all it's zero-valued entries indexed by $\mathcal{S}^{\rm c}$ removed).  Let ${\bf \tilde{f}'}_{\mathcal{S}}, {\bf z}_{\rm min} \in \mathbbm{R}^{|\mathcal{S}|}$ be defined as in line 2 of Algorithm~\ref{alg:FSL}.  We have that
\begin{equation}
\begin{array}{llll} 
\left \|{\bf \tilde{f}'}_{\mathcal{S}} - {\bf \tilde{f}''}_{\mathcal{S}} \right \|_2 &\leq~ \frac{\left \| \widetilde{R}_{\mathcal{S}} {\bf \tilde{f}'}_{\mathcal{S}} - \widetilde{R}_{\mathcal{S}}  {\bf \tilde{f}''}_{\mathcal{S}} \right \|_2}{\sqrt{1-\epsilon}} & \hspace{.5in} & \textrm{Since \eqref{thm:randSampeq} holds}\\
&\leq~ \frac{\left \| \widetilde{R}_{\mathcal{S}} {\bf z}_{\rm min} - \widetilde{R}_{\mathcal{S}}  {\bf \tilde{f}''}_{\mathcal{S}} \right \|_2 + ~\delta}{\sqrt{1-\epsilon}} & \hspace{.5in} & \textrm{By CG discussion \eqref{equ:ApproxCGsolve} in \S\ref{sec:Alg}}\\
&\leq ~\frac{1}{\sqrt{1-\epsilon}} \left( \left \| \widetilde{R}_{\mathcal{S}} {\bf z}_{\rm min} - {\bf y} \right \|_2 + \left \| \widetilde{R} {\bf \tilde{f}}_{\mathcal{S}^{\rm c}} \right \|_2 + ~\delta \right) &  &\textrm{Since }{\bf y} = \widetilde{R} \left( {\bf \tilde{f}}_{\mathcal{S}} + {\bf \tilde{f}}_{\mathcal{S}^{\rm c}} \right)\\
&\leq  ~\frac{2}{\sqrt{1-\epsilon}} \left \| \widetilde{R} {\bf \tilde{f}}_{\mathcal{S}^{\rm c}} \right \|_2 + \frac{\delta}{\sqrt{1-\epsilon}} &  &\textrm{By the definition of } {\bf z}_{\rm min}\\
&\leq~\frac{2 \sqrt{1+\epsilon}}{\sqrt{1-\epsilon}} \left( \left \| {\bf \tilde{f}}_{\mathcal{S}^{\rm c}} \right \|_2 + \frac{\left \| {\bf \tilde{f}}_{\mathcal{S}^{\rm c}} \right \|_1}{\sqrt{s}} \right) + \frac{\delta}{\sqrt{1-\epsilon}} & &\textrm{Using Exercise 6.6 in \cite{HolgerBook}}
\end{array}
\nonumber .
\end{equation}
One can now (implicitly) form ${\bf z} = {\bf z}_{\mathcal{S}}$ from ${\bf \tilde{f}'}_{\mathcal{S}}$ in the obvious fashion.  Doing so we learn that
\begin{align}
\left \| {\bf \tilde{f}} - {\bf z}  \right \|_2 & \leq  \left \| {\bf z}_{\mathcal{S}} - {\bf \tilde{f}}_{\mathcal{S}} \right \|_2 + \left \| {\bf \tilde{f}}_{\mathcal{S}^{\rm c}}  \right \|_2 \nonumber \\
& \leq \left( 1 + \frac{2 \sqrt{1+\epsilon}}{\sqrt{1-\epsilon}}  \right) \left \| {\bf \tilde{f}}_{\mathcal{S}^{\rm c}} \right \|_2 + \frac{2 \sqrt{1+\epsilon}}{\sqrt{1-\epsilon}} \cdot \frac{\left \| {\bf \tilde{f}}_{\mathcal{S}^{\rm c}} \right \|_1}{\sqrt{s}}  + \frac{\delta}{\sqrt{1-\epsilon}} \nonumber .
\end{align}
The total runtime complexity follows from the discussion regarding line 2 of Algorithm~\ref{alg:FSL} in Section~\ref{sec:Alg}.  \qed\\

We are now prepared to test a particular version of Algorithm~\ref{alg:FSL} numerically.  As we shall see, the experiments demonstrate (as a proof of concept) that SFTs can be used to build stable sublinear-time algorithms capable of rapidly computing Legendre coefficients whenever they exhibit compressibility.

\section{Empirical Evaluation}
\label{sec:Eval}

We now present representative results demonstrating the numerical robustness and efficiency of the proposed SFT-based Legendre method.  For the experiments in this section we used FFTW 3.3.4 to implement Iserles' scheme \cite{iserles2011fast} for the purposes of comparison (see \eqref{equ:IserlesOriginal} in \S\ref{sec:IserlesMap}).  FFTW3 \cite{FFTW3ref} is a highly  efficient implementation of the ``standard'' FFT algorithm -- it has been systematically optimized over the course of the last two decades and remains one of the fastest freely available FFT implementations available today.  The parameter $r$ in \eqref{def:Fmap} was chosen to be $1 - 10^{-8}$ for all experiments.  The qualitative behavior for other values of $r$ sufficiently close to $1$ is similar.  The number of terms, $M$, in \eqref{equ:IserlesOriginal} was varied differently in each experiment, as indicated below.

Algorithm~\ref{alg:SFTLegendre} was implemented using AAFFT \cite{AAFFT1exp} as the sparse Fourier transform, followed by a conjugate gradient code\footnote{The conjugate gradient code is available at \url{http://people.sc.fsu.edu/~jburkardt/cpp_src/cg/cg.html}.} in order to compute the necessary Legendre coefficients.  Although less optimized than some other SFT implementations, AAFFT's code is well documented, readable, and easy to modify.  All code used to perform the experiments below is freely available.\footnote{All code is available at \url{http://www.math.msu.edu/~markiwen/Code.html}.}

\begin{figure}[hbtp]
        \centering
        \begin{subfigure}[b]{0.495\textwidth}
                \includegraphics[width=1\textwidth]{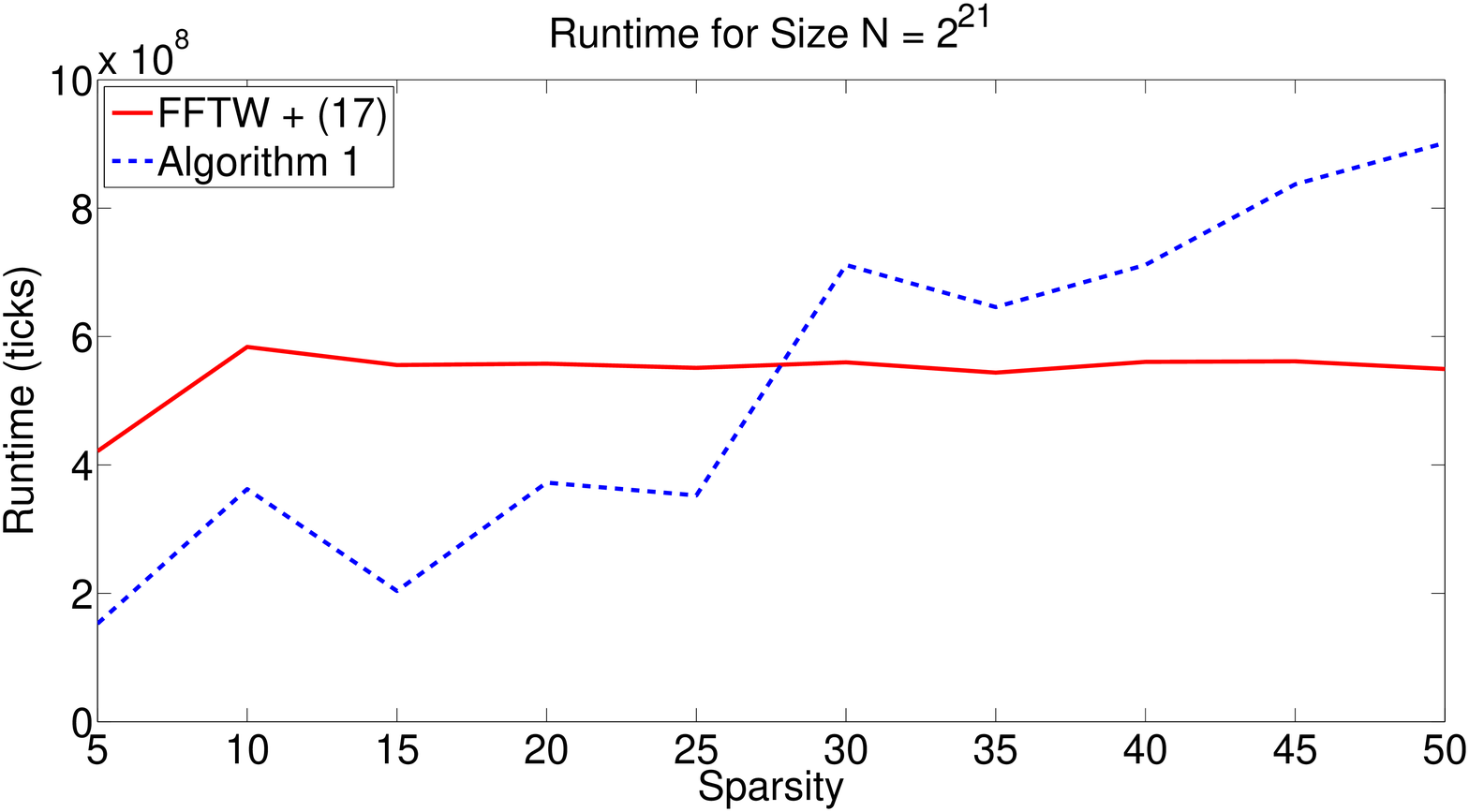}
                \caption{Runtime with $N = 2^{21}$ as sparsity, $s$, varies.}
                \label{fig:runtimeSvaries}
        \end{subfigure}
        \hfill
        \begin{subfigure}[b]{0.495\textwidth}
                \includegraphics[width=1\textwidth]{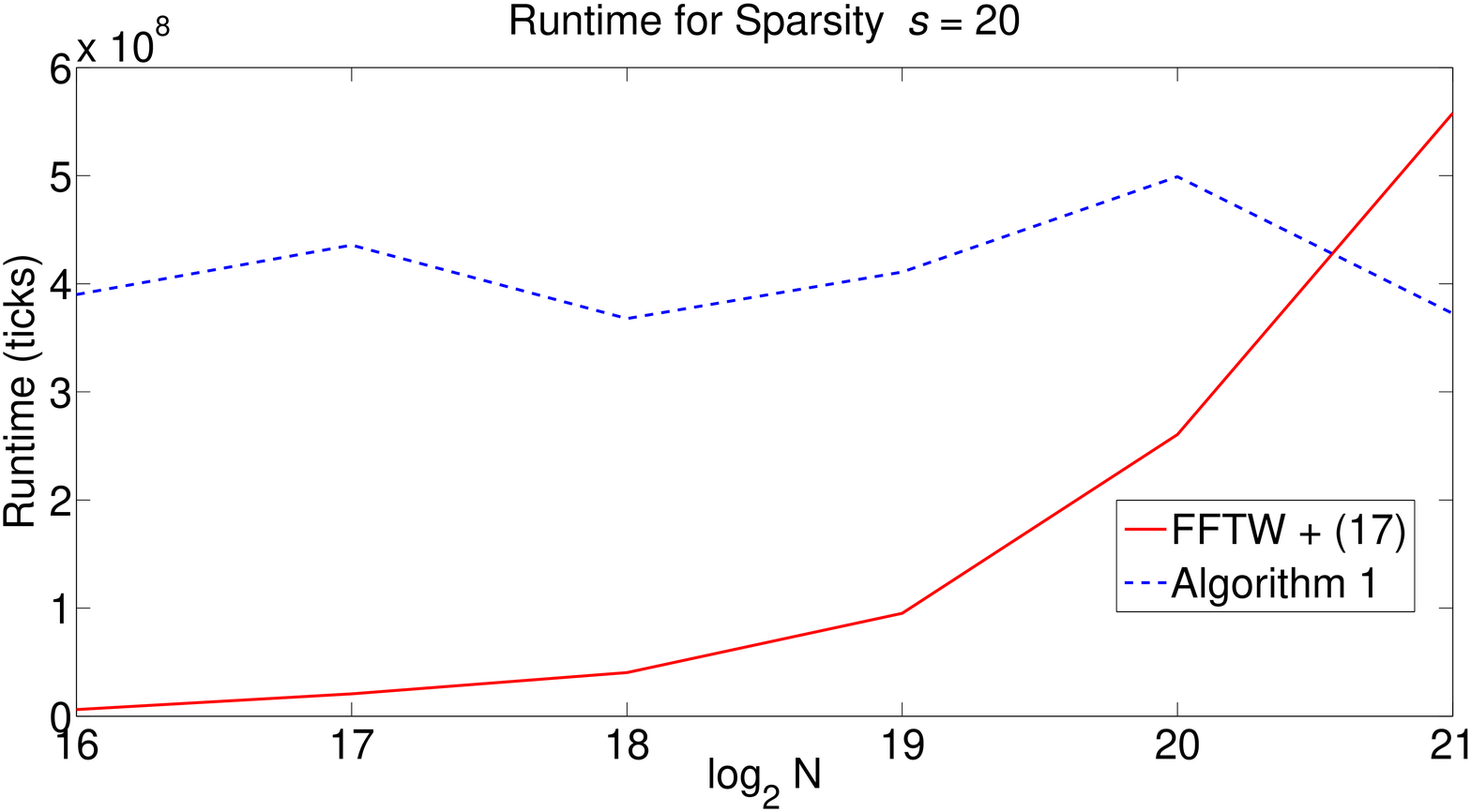}
                \caption{Runtime with $s = 20$ as $N$ varies.}
                \label{fig:runtimeNvaries}
        \end{subfigure}
        \caption{Runtime comparison between Algorithm~\ref{alg:SFTLegendre} implemented with AAFFT as the SFT, and the approach from \cite{iserles2011fast} implemented with FFTW as the FFT.}
        \label{fig:efficiency}
\end{figure}

Every data point in the first two figures below is the result of 100 trials performed on 100 different randomly generated polynomials,
\begin{equation}
f(x) = \sum_{m \in S,~|S| = s} a_m L_m(x),
\label{equ:TrialSig}
\end{equation}
where $S \subset [N]$ contains $s$ entries independently chosen uniformly at random from $[N]$, and each $a_m \in \{ -1, 1\}$ is independently chosen to be $\pm 1$ with probability $1/2$.  Figure~\ref{fig:efficiency} reports runtime results averaged over the 100 independent trials at each data point for Algorithm~\ref{alg:SFTLegendre} versus the method outlined in \S\ref{sec:IserlesMap}.  All runtimes are reported in tick counts using the ``cycle.h'' header included with the FFTW 3.3.4 library.  Tick counts correspond closely to the number of CPU cycles used during the execution of a program and, therefore, accurately reflect each implementation's comparative computational complexity.  
As one can see, Algorithm~\ref{alg:SFTLegendre} is faster than the method outlined in \S\ref{sec:IserlesMap} when $s \ll N$.\footnote{Using a faster SFT implementation would doubtlessly produce faster results for Algorithm~\ref{alg:SFTLegendre} -- AAFFT is the computational bottleneck here.}

\begin{figure}[hbtp]
        \centering
        \begin{subfigure}[b]{0.495\textwidth}
                \includegraphics[width=1\textwidth]{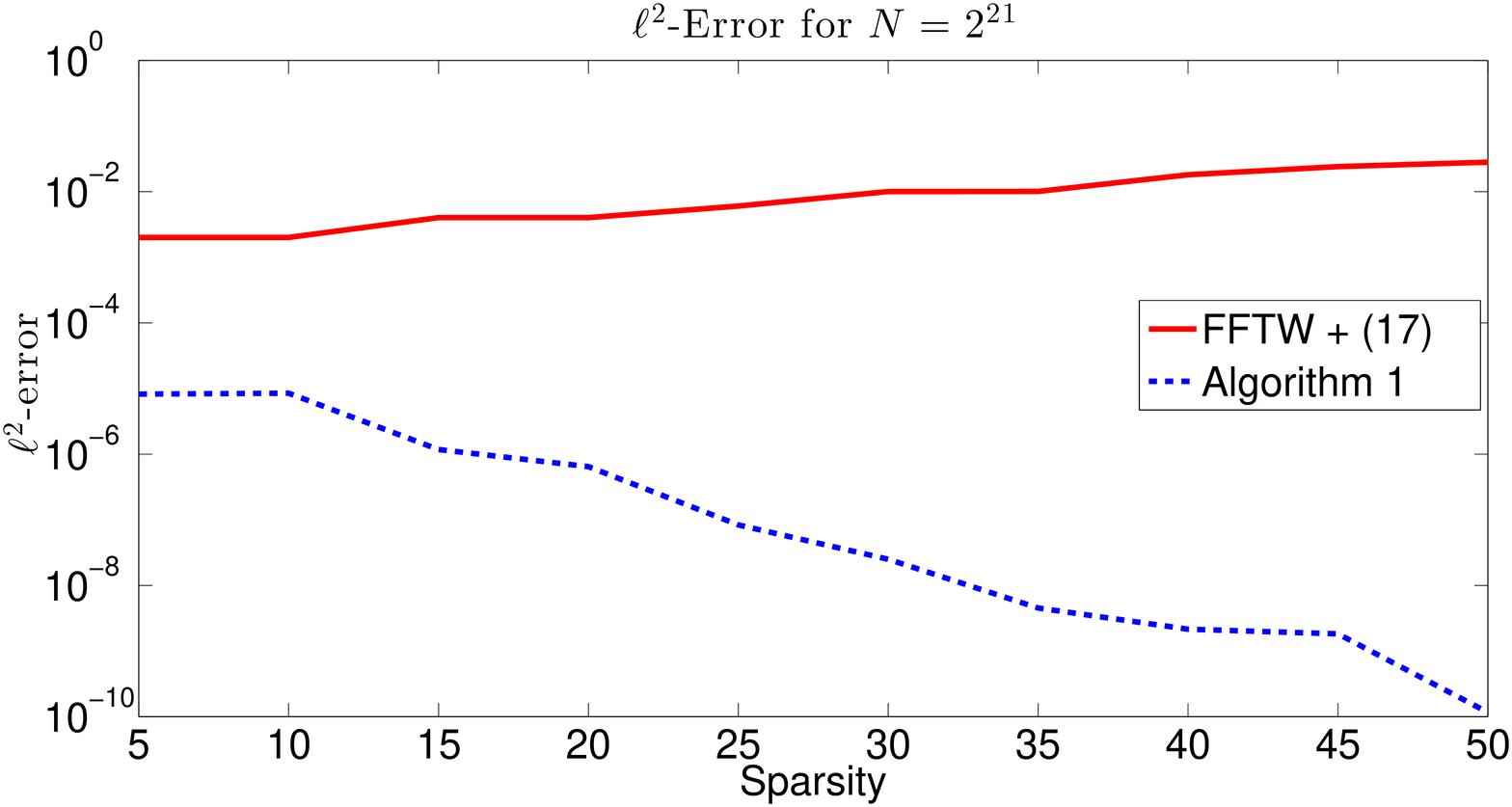}
                \caption{Average $\ell^2$ error for $N = 2^{21}$ runtime experiments}
                \label{fig:errorSvaries}
        \end{subfigure}
        \hfill
        \begin{subfigure}[b]{0.495\textwidth}
                \includegraphics[width=1\textwidth]{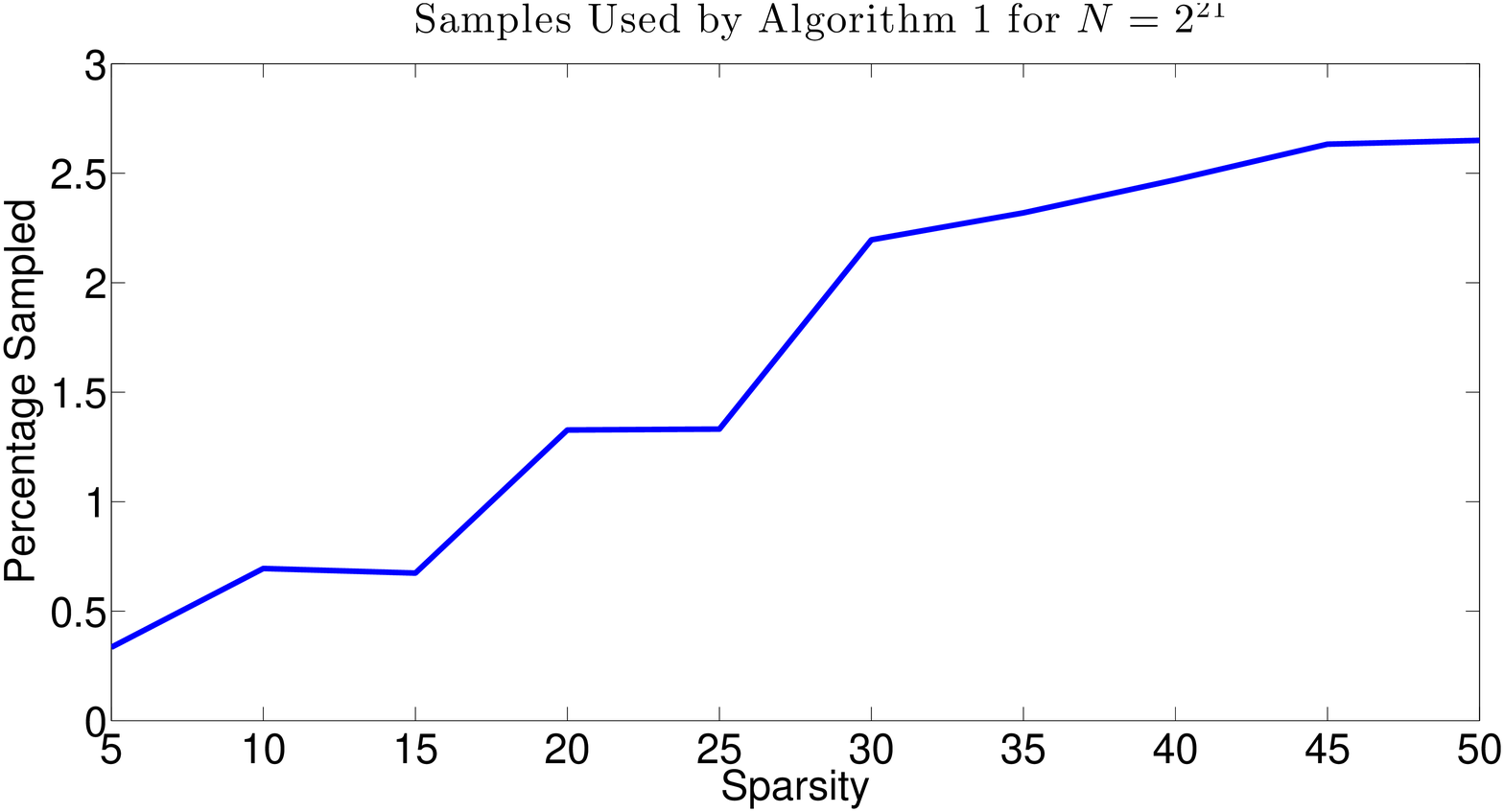}
                \caption{Algorithm~\ref{alg:SFTLegendre} function evaluations for $N = 2^{21}$}
                \label{fig:AAFFTsamps}
        \end{subfigure}
        \caption{Additional information regarding the runtime experiments in Figure~\ref{fig:runtimeSvaries}.  The number of terms, $M$, used in \eqref{equ:IserlesOriginal} for Figure~\ref{fig:errorSvaries} was chosen separately from the range $[1, 25]$ for each of the 100 trials in order to give the lowest possible error.  The best value of $M$ was usually $1$, however.  Choosing $M$ to be much larger did not usually help.  Figure~\ref{fig:AAFFTsamps} plots the average number of evaluations of each trial polynomial $f$, as a percentage of $N$, that Algorithm~\ref{alg:SFTLegendre} used in order to produce the results plotted in Figures~\ref{fig:errorSvaries} and~\ref{fig:runtimeSvaries}.}
        \label{fig:error1}
\end{figure}

AAFFT is an implementation of a randomized Fourier method with a user-tunable probability of failure (see Theorem~\ref{thm:MultDimRecov} for a similar probabilistic SFT recovery guarantee).  The runtime results for Algorithm~\ref{alg:SFTLegendre} in Figure~\ref{fig:efficiency} were produced with settings for AAFFT which guaranteed Algorithm~\ref{alg:SFTLegendre} to have an $\ell^2$-error of size $< 10^{-5}$ on the $s$ true Legendre coefficients, $a_m$ with $m \in S$, for more than $70$ of the $100$ trials used to generate each data point.  More detailed information concerning approximation errors for these experiments is reported in Figure~\ref{fig:errorSvaries}.  Although highly accurate for small values of $N$ (see \cite{iserles2011fast}), the method outlined in \S\ref{sec:IserlesMap} has a relatively low accuracy for the large degree polynomials considered herein, achieving only two or three digits of accuracy per Legendre coefficient on average.  Figure~\ref{fig:errorSvaries} compares its average $\ell^2$-error on the true Legendre coefficients, $a_m$ with $m \in S$, over all 100 trials at each data point against Algorithm~\ref{alg:SFTLegendre}'s average $\ell^2$-error over the at least 70 trials at each data point for which it correctly identified a superset of $S$ from \eqref{equ:TrialSig}.  As one can see, Algorithm~\ref{alg:SFTLegendre} is generally more accurate \textit{when it manages to identify $S$}.  The error graphs for the other values of $N$ considered herein were similar.  For example, the method outlined in \S\ref{sec:IserlesMap} had an average $\ell^2$-error that was always less than $0.05$ for the experiments reported in Figure~\ref{fig:runtimeNvaries} at each $N$, while Algorithm~\ref{alg:SFTLegendre}'s average $\ell^2$-error was always $< 10^{-5}$ for these experiments whenever it found $S$ (for more than $70\%$ of the trials for each data point).

\begin{figure}[hbtp]
\centering
\includegraphics[width=1\textwidth]{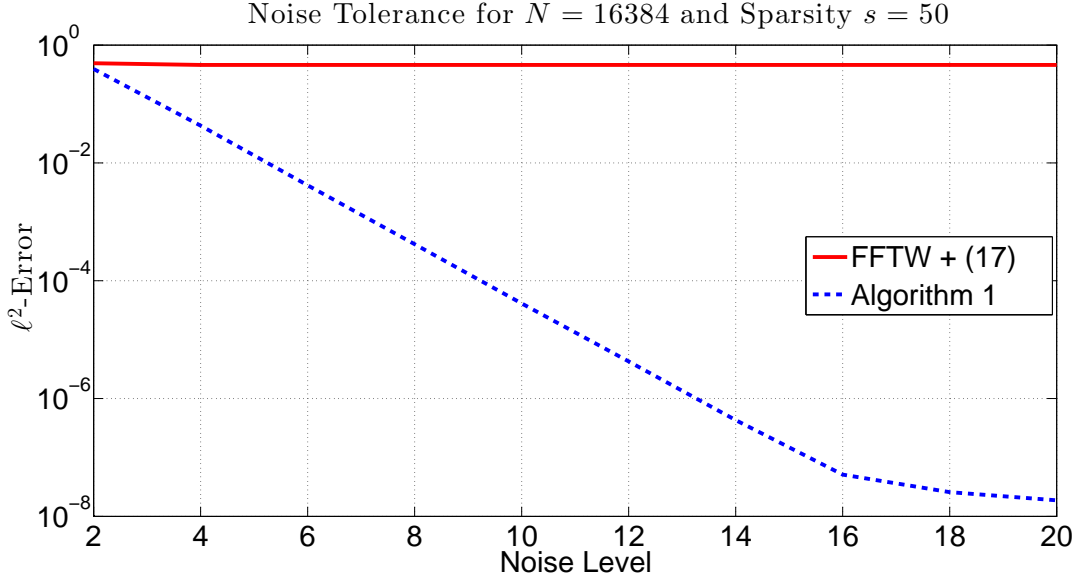}
\caption{Robustness to approximate sparsity.  The horizontal axis varies over ten signal-to-noise values given by $\log_{10} \left( \frac{s}{\sum_{m \in [N] \setminus S} b^2_m} \right) ~=~ \log_{10} \left( \frac{50}{\sum_{m \in [N] \setminus S} b^2_m} \right)$.  The vertical axis graphs the average $\ell^2$-error of each method.}
\label{fig:robust}
\end{figure}

Figure~\ref{fig:robust} reports the results of some additional experiments on numerical accuracy, stability, and robustness to noise.  For these experiments both the maximum degree, $N$, and the sparsity, $s$, of the trial functions were fixed.  In addition, the trial functions were modified so that each one was of the form
\begin{equation}
f(x) = \sum_{m \in S,~|S|=s} a_m L_m(x) + \sum_{m \in [N] \setminus S} b_m L_m(x),
\label{equ:TrialSig2}
\end{equation}
where $(i)$ $S \subset [N]$ contains $s$ entires independently chosen uniformly at random from $[N]$, $(ii)$ each $a_m \in \{ -1, 1\}$ is independently chosen to be $\pm 1$ with probability $1/2$, and $(iii)$ each $b_m$ is an i.i.d mean 0 Gaussian random number generated numerically via the Box-Muller method.  As above, each data point in Figure~\ref{fig:robust} is the result of 100 trials performed on 100 independently generated polynomials of this form \eqref{equ:TrialSig2}.

For the experiments in Figure~\ref{fig:robust}, Algorithm~\ref{alg:SFTLegendre} had its parameters set so that it would compute all Legendre coefficients to at least $8$ digits of precision for at least $90\%$ of trials in the noiseless setting (i.e., when \eqref{equ:TrialSig} holds).  The log signal-to-noise ratio,
\begin{equation}
\log_{10} \left( \frac{\sum_{m \in S} a^2_m }{\sum_{m \in [N] \setminus S} b^2_m} \right) ~=~ \log_{10} \left( \frac{s}{\sum_{m \in [N] \setminus S} b^2_m} \right),
\end{equation}
was then varied by renormalizing the i.i.d. $b_m$'s generated for each trial.  The resulting log signal-to-noise ratios appear on the horizontal axis of Figure~\ref{fig:robust}.  The vertical axis plots the average $\ell^2$-error on the true Legendre coefficients, $a_m$ with $m \in S$, over all trials at each noise level of the method in \S\ref{sec:IserlesMap} against Algorithm~\ref{alg:SFTLegendre}'s average $\ell^2$-error over the at least 90 trials at each noise level for which it correctly identified a superset of $S \subset [N]$ from \eqref{equ:TrialSig2}.  The number of terms, $M$, used in \eqref{equ:IserlesOriginal} for the \S\ref{sec:IserlesMap} method was again chosen separately for each of the 100 trials in order to give the lowest possible error, except that here $M \in [1, 45]$.  The best value of $M$ was, as before, usually $1$, however.  Increasing $M$ still did not usually help decrease the error.  The qualitative behavior for other values of $s$ and $N$ was similar.  As above, we conclude that Algorithm~\ref{alg:SFTLegendre} is generally more accurate than the original  \S\ref{sec:IserlesMap} method provided that it correctly identifies (a superset of) the support set $S$.

\section{Conclusion}
\label{sec:Conc}

In this paper we have demonstrated that SFT techniques can be used to help rapidly approximate functions with sparse Legendre coefficient expansions.  Together with the problems already relegated to future work, we believe that it would be interesting to consider extending these methods to functions which exhibit sparsity in other types of Gegenbauer polynomial expansions.  Given the level of success in both the Chebyshev and Legendre settings it seems likely that SFTs can be utilized to good effect more generally.

\section*{Acknowledgements}

The authors would like to thank Aditya Viswanathan for many helpful suggestions and discussions during the writing of this paper.

\bibliographystyle{abbrv}
\bibliography{Legendre}

\appendix

\section{Proof of Lemma~\ref{lem:InverseMatrix}}
\label{app:ProofLemMatrixInv}

Assuming that $f:  [-1, 1] \rightarrow \mathbbm{R}$ is a Legendre polynomial of degree $N$, we can begin to compute $\widehat{f}_r$.  This will, in turn, reveal the entries of $\tilde{G}^{-1}_r$ from \eqref{equ:IserlesNew}.  Recall that
\begin{align}
f_r( x ) &= \left( 1 - r^2 \mathbbm{e}^{2 \mathbbm{i} x} \right) f\left( \frac{1}{2}\left( r^{-1} \mathbbm{e}^{- \mathbbm{i} x} + r \mathbbm{e}^{\mathbbm{i} x} \right) \right) \nonumber \\
& = \left( 1 - r^2 \mathbbm{e}^{2 \mathbbm{i} x} \right) \left[ \sum^{\left \lfloor \frac{N}{2} \right\rfloor}_{j=0} \tilde{f}(2j) \cdot L_{2j}\left( \frac{ r^{-1} \mathbbm{e}^{- \mathbbm{i} x} + r \mathbbm{e}^{\mathbbm{i} x}}{2}\right) + \sum^{\left \lceil \frac{N}{2} - 1\right \rceil}_{j=0} \tilde{f}(2j+1) \cdot L_{2j+1}\left( \frac{ r^{-1} \mathbbm{e}^{- \mathbbm{i} x} + r \mathbbm{e}^{\mathbbm{i} x}}{2} \right)\right] .  \nonumber
\end{align}
Setting $z := \frac{ r^{-1} \mathbbm{e}^{- \mathbbm{i} x} + r \mathbbm{e}^{\mathbbm{i} x}}{2}$ and expanding the Legendre polynomials in terms of the canonical polynomial basis (see, e.g., \cite{zhang1996computation}) we obtain
\begin{align}
f_r( x ) = \left( 1 - r^2 \mathbbm{e}^{2 \mathbbm{i} x} \right) &\left[ \sum^{\left \lfloor \frac{N}{2} \right\rfloor}_{j=0} \tilde{f}(2j)
\sum^{j}_{k=0} (-1)^k \frac{(4j-2k)!}{4^j k! (2j - k)! (2j - 2k)!} z^{2(j-k)} \right.\\
&\left.+ \sum^{\left \lceil \frac{N}{2} - 1\right \rceil}_{j=0} \tilde{f}(2j+1) \sum^j_{k=0} (-1)^k \frac{(4j-2k+2)!}{2 \cdot 4^j k! (2j - k +1)! (2j - 2k +1)!} z^{2(j-k)+1} \right] .  \nonumber
\end{align}
Note that even powers of $z$ will only ever contain even powers of both $r^{-1} \mathbbm{e}^{- \mathbbm{i} x}$ and $r \mathbbm{e}^{\mathbbm{i}x}$.  Similarly, odd powers of $z$ will only ever contain odd powers of $r^{-1} \mathbbm{e}^{- \mathbbm{i} x}$ and $r \mathbbm{e}^{\mathbbm{i}x}$.  Thus, it we can see that
$$\widehat{f}_r (\omega) = \left\{ \begin{array}{ll} \widehat{f}_e (\omega) & {\rm if}~ \omega \equiv 0~{\rm mod}~2 \\ \widehat{f}_o (\omega) & {\rm if}~ \omega \equiv 1~{\rm mod}~2 \end{array} \right.$$
where
\begin{align}
f_e( x ) &= \left( 1 - r^2 \mathbbm{e}^{2 \mathbbm{i} x} \right) \left[ \sum^{\left \lfloor \frac{N}{2} \right\rfloor}_{j=0} \tilde{f}(2j)
\sum^{j}_{k=0} (-1)^k \frac{(4j-2k)!}{4^j k! (2j - k)! (2j - 2k)!} \left( \frac{ r^{-1} \mathbbm{e}^{- \mathbbm{i} x} + r \mathbbm{e}^{\mathbbm{i} x}}{2}\right)^{2(j-k)} \right] \nonumber \\
&= \left( 1 - r^2 \mathbbm{e}^{2 \mathbbm{i} x} \right) \left[ \sum^{\left \lfloor \frac{N}{2} \right\rfloor}_{j=0} \tilde{f}(2j) \frac{(-1)^j}{4^j}
\sum^{j}_{k=0} (-1)^k \frac{(2j+2k)!}{(j-k)! (j + k)! (2k)!} \left( \frac{ r^{-1} \mathbbm{e}^{- \mathbbm{i} x} + r \mathbbm{e}^{\mathbbm{i} x}}{2}\right)^{2k} \right]
\label{equ:EVENcoefs}
\end{align}
and
\begin{align}
f_o( x ) = &\left( 1 - r^2 \mathbbm{e}^{2 \mathbbm{i} x} \right) \cdot \nonumber \\ &\left[ \sum^{\left \lceil \frac{N}{2} - 1\right \rceil}_{j=0} \tilde{f}(2j+1) \sum^j_{k=0} (-1)^k \frac{(4j-2k+2)!}{2 \cdot 4^j k! (2j - k +1)! (2j - 2k +1)!} \left( \frac{ r^{-1} \mathbbm{e}^{- \mathbbm{i} x} + r \mathbbm{e}^{\mathbbm{i} x}}{2}\right)^{2(j-k)+1} \right] \nonumber \\
= &\left( 1 - r^2 \mathbbm{e}^{2 \mathbbm{i} x} \right) \cdot \nonumber \\ &\left[ \sum^{\left \lceil \frac{N}{2} - 1\right \rceil}_{j=0} \tilde{f}(2j+1) \frac{(-1)^j}{2 \cdot 4^j} \sum^j_{k=0} (-1)^k \frac{(2j+2k+2)!}{(j-k)! (j + k +1)! (2k +1)!} \left( \frac{ r^{-1} \mathbbm{e}^{- \mathbbm{i} x} + r \mathbbm{e}^{\mathbbm{i} x}}{2}\right)^{2k+1} \right].
\label{equ:ODDcoefs}
\end{align}

\subsection{The Even Fourier Coefficients}

Expanding $f_e(x)$ from \eqref{equ:EVENcoefs} using the Binomial Theorem, reindexing, and then changing the order of summation, we obtain
\begin{align}
f_e( x ) &= \left( 1 - r^2 \mathbbm{e}^{2 \mathbbm{i} x} \right) \left[ \sum^{\left \lfloor \frac{N}{2} \right\rfloor}_{j=0} \tilde{f}(2j) \frac{(-1)^j}{4^j}
\sum^{j}_{k=0} \frac{(-1)^k}{4^k} \frac{(2j+2k)!}{(j-k)! (j + k)! (2k)!} \sum^k_{l = -k}  {2k \choose k+l} \left( r \mathbbm{e}^{\mathbbm{i} x} \right)^{2l} \right] \nonumber \\
&=  \left( 1 - r^2 \mathbbm{e}^{2 \mathbbm{i} x} \right) \left[ \sum^{\left \lfloor \frac{N}{2} \right\rfloor}_{l = -\left \lfloor \frac{N}{2} \right\rfloor} \left( r \mathbbm{e}^{\mathbbm{i} x} \right)^{2l} \sum^{\left \lfloor \frac{N}{2} \right\rfloor}_{j=|l|} \tilde{f}(2j) \frac{(-1)^j}{4^j} \sum^{j}_{k=|l|} \frac{(-1)^k}{4^k} \frac{(2j+2k)!}{(j-k)! (j + k)! (2k)!} {2k \choose k+l} \right]. \nonumber
\end{align}
Multiplying the $\left( 1 - r^2 \mathbbm{e}^{2 \mathbbm{i} x} \right)$ factor through the sum above and recombining terms now yields
\begin{align}
f_e(x)=&\left[ \left(r \mathbbm{e}^{\mathbbm{i} x}\right)^{-2 \left\lfloor \frac{N}{2} \right\rfloor}- \left( r \mathbbm{e}^{\mathbbm{i} x} \right)^{2 \left \lfloor \frac{N}{2} \right\rfloor + 2} \right] {4 \left \lfloor \frac{N}{2} \right\rfloor \choose 2\left \lfloor \frac{N}{2} \right\rfloor} \frac{\tilde{f} \left(2 \left \lfloor N / 2 \right\rfloor \right)}{16^{\left \lfloor N / 2 \right\rfloor}} \nonumber \\
&+ \sum^{\left \lfloor \frac{N}{2} \right\rfloor}_{l = -\left \lfloor \frac{N}{2} \right\rfloor+1} \left( r \mathbbm{e}^{\mathbbm{i} x} \right)^{2l}\left \{ \sum^{\left \lfloor \frac{N}{2} \right\rfloor}_{j=|l|} \tilde{f}(2j) \frac{(-1)^j}{4^j} \sum^{j}_{k=|l|} \frac{(-1)^k}{4^k} \frac{(2j+2k)!}{(j-k)! (j + k)! (2k)!} {2k \choose k+l} \right. \label{equ:EvenLastTerm} \\ & \hspace{1.5in} \left. - \sum^{\left \lfloor \frac{N}{2} \right\rfloor}_{j= |l-1|}  \tilde{f}(2j) \frac{(-1)^{j}}{4^{j}} \sum^{j}_{k=|l-1|} \frac{(-1)^k}{4^k} \frac{(2j+2k)!}{(j-k)! (j + k)! (2k)!} {2k \choose k+l} \frac{k+l}{k+1-l}\right \}\nonumber
\end{align}

Recalling from \eqref{equ:IserlesNew} that we are primarily concerned with $l = 0, -1, \dots, -\left \lfloor \frac{N}{2} \right\rfloor$ in the expression above, we can now recombine terms in the bracketed sum to obtain the relevant even Fourier series coefficients of $f_r$ from $f_e$ (please note that, in fact, that the last line of our calculation above does not technically hold as written unless $l \leq 0$!).  Doing so, we learn that
\begin{align}
\widehat{f}_r \left( 2 l \right) = & r^{2l} \left\{ \tilde{f}(2 |l|) \frac{1}{16^{|l|}} {4 |l| \choose 2 |l|} +\sum^{\left \lfloor \frac{N}{2} \right\rfloor }_{j=|l|+1} \tilde{f}(2 j)  \frac{(-1)^j}{4^j} \left[ \frac{(-1)^{|l|}}{4^{|l|}}\frac{ (2j+2|l|)!}{(j-|l|)!(j+|l|)!(2|l|)!} \right. \right.\nonumber \\
&+ \left. \left. \sum^j_{k=|l|+1}  \frac{(-1)^k}{4^k} \frac{(2j+2k)!}{(j-k)! (j + k)! (2k)!} {2k \choose k+l} \frac{1-2l}{k-l+1} \right]  \right\}
\end{align}
for $l = 0, -1, \dots, -\left \lfloor \frac{N}{2} \right\rfloor+1$.  This combined with \eqref{equ:EvenLastTerm} gives all entries in the even rows of $\tilde{G}^{-1}_r$.

\subsection{The Odd Fourier Coefficients}

Expanding $f_o(x)$ from \eqref{equ:ODDcoefs} using the Binomial Theorem, reindexing, and then changing the order of summation, we obtain
\begin{align}
f_o( x ) =& \left( 1 - r^2 \mathbbm{e}^{2 \mathbbm{i} x} \right) \cdot \nonumber \\ &\left[ \sum^{\left \lceil \frac{N}{2} - 1\right \rceil}_{j=0} \tilde{f}(2j+1) \frac{(-1)^j}{2 \cdot 4^j} \sum^j_{k=0} \frac{(-1)^k}{2 \cdot 4^k} \frac{(2j+2k+2)!}{(j-k)! (j + k +1)! (2k +1)!} \sum^{k+1}_{l = -k}  {2k+1 \choose k+l} \left( r \mathbbm{e}^{\mathbbm{i} x} \right)^{2l-1}  \right]\nonumber \\
=& \left( 1 - r^2 \mathbbm{e}^{2 \mathbbm{i} x} \right) \cdot \left[ \sum^{\left \lceil \frac{N}{2} - 1\right \rceil}_{l=0} C_l \left( r \mathbbm{e}^{\mathbbm{i} x} \right)^{2l+1} + \right. \nonumber \\
&\left. \sum^{0}_{l = -\left \lceil \frac{N}{2}-1 \right\rceil} \left( r \mathbbm{e}^{\mathbbm{i} x} \right)^{2l-1} \sum^{\left \lceil \frac{N}{2}-1 \right\rceil}_{j=|l|} \tilde{f}(2j+1) \frac{(-1)^j}{2\cdot4^j} \sum^{j}_{k=|l|} \frac{(-1)^k}{2\cdot 4^k} \frac{(2j+2k+2)!}{(j-k)! (j + k+1)! (2k+1)!} {2k+1 \choose k+l} \right]. \nonumber
\end{align}
for some $C_0, \dots, C_{\left \lceil \frac{N}{2} - 1\right \rceil} \in \mathbbm{R}$ with which we need not concern ourselves at the moment.  Multiplying the $\left( 1 - r^2 \mathbbm{e}^{2 \mathbbm{i} x} \right)$ factor through the sum above and recombining terms now yields
\begin{align}
f_o( x ) =& \left( 1 - r^2 \mathbbm{e}^{2 \mathbbm{i} x} \right)\sum^{\left \lceil \frac{N}{2} - 1\right \rceil}_{l=0} C_l \left( r \mathbbm{e}^{\mathbbm{i} x} \right)^{2l+1} +  \nonumber \label{equ:OddLastTerm}\\
& 
\left(r \mathbbm{e}^{\mathbbm{i} x}\right)^{-2 \left\lceil \frac{N}{2}-1 \right\rceil-1}
{4 \left \lceil \frac{N}{2}-1 \right\rceil +2 \choose 2\left \lceil \frac{N}{2}-1 \right\rceil +1} \frac{\tilde{f} \left(2 \left \lceil N / 2-1 \right\rceil +1 \right)}{4\cdot 16^{\left \lceil N / 2 -1\right\rfloor}} \\
+&\sum^{0}_{l = -\left \lceil \frac{N}{2}-1 \right\rceil+1} \left( r \mathbbm{e}^{\mathbbm{i} x} \right)^{2l-1} \left\{ \sum^{\left \lceil \frac{N}{2}-1 \right\rceil}_{j=|l|} \tilde{f}(2j+1) \frac{(-1)^j}{2\cdot4^j}\sum^{j}_{k=|l|} \frac{(-1)^k}{2\cdot 4^k} \frac{(2j+2k+2)!}{(j-k)! (j + k+1)! (2k+1)!} {2k+1 \choose k+l} \right. \nonumber \\
& \hspace{0.5in} \left. - \sum^{\left \lceil \frac{N}{2}-1 \right\rceil}_{j= |l-1|}  \tilde{f}(2j+1) \frac{(-1)^{j}}{2\cdot 4^{j}} \sum^{j}_{k=|l-1|} \frac{(-1)^k}{2\cdot 4^k} \frac{(2j+2k+2)!}{(j-k)! (j + k+1)! (2k+1)!} {2k+1 \choose k+l} \frac{k+l}{k+2-l}\right\}.\nonumber
\end{align}

We can now recombine terms in the bracketed sum to obtain the relevant odd Fourier series coefficients of $f_r$ from $f_o$.  Doing so, we learn that
\begin{align}
\widehat{f}_r \left( 2 l -1\right) = & \,\, r^{2l-1} \left\{ \tilde{f}(2 |l|+1) \frac{1}{4\cdot 16^{|l|}} {4 |l|+2 \choose 2 |l|+1}  \right. \nonumber \\
& +\sum^{\left \lceil \frac{N}{2}-1 \right\rceil }_{j=|l|+1} \tilde{f}(2 j+1)  \frac{(-1)^j}{2\cdot 4^j} \left[ \frac{(-1)^{|l|}}{2\cdot 4^{|l|}}\frac{ (2j+2|l|+2)!}{(j-|l|)!(j+|l|+1)!(2|l|+1)!} \right. \nonumber \\
&\left. \left. + \sum^j_{k=|l|+1}  \frac{(-1)^k}{2\cdot 4^k} \frac{(2j+2k+2)!}{(j-k)! (j + k+1)! (2k+1)!} {2k+1 \choose k+l} \frac{2-2l}{k-l+2} \right]  \right\}
\end{align}
for $l = 0, -1, \dots, -\left \lceil \frac{N}{2}-1 \right\rceil+1$.  This combined with \eqref{equ:OddLastTerm} gives all entries in the odd rows of $\tilde{G}^{-1}_r$.

\section{Proof of Corollary~\ref{cor:CompactInverseMatrix}}
\label{app:ProofLemCompactMatrixInv}

We will once again consider the even and odd rows separately.  

\subsection{The Even Rows}
From \eqref{EvenrowsofInverse} in Lemma ~\ref{lem:InverseMatrix}, the nonzero entries in the even rows of the inverse matrix $\tilde{G}_{r}^{-1}$ for $i \leq j \leq \left \lfloor \frac{N}{2} \right\rfloor$ can be rewritten as 
\begin{align}
\left(\tilde{G}_{r}^{-1}\right)_{2i, 2j} 
= & \left( -\frac{1}{4}\right)^{i+j} r^{-2i}(1+2i)\sum_{k=0}^{j-i}\left(-\frac{1}{4}\right)^{k} \frac{(2j+2i+2k)!}{\left((j-i)-k\right)!(j+i+k)!(2i+1+k)!k!}. \nonumber
\end{align}
For every $i, j, k\ge 0$ and $m\in \mathbb{Z}_{+}$, 
\begin{equation}\label{identities}
\begin{array}{ll}
& (-m)_{k}=\frac{(-1)^{k}m!}{(m-k)!}, \quad (m+k)!=(m)_{k+1}(m-1)!, \quad {\rm and} \\ 
&\left(2(i+j)+2k\right)! =4^k\left(2(j+i)\right)!(i+j+1)_{k}\left(i+j+\frac{1}{2}\right)_{k}. 
\end{array}
\end{equation}
Thus, we deduce that
\begin{align}
\left(\tilde{G}_{r}^{-1}\right)_{2i, 2j} =& \left(-\frac{1}{4}\right)^{i+j}\frac{r^{-2i}(1+2i)}{(j-i)!}\sum_{k=0}^{j-i}\frac{\left(-(j-i)\right)_{k}(2j+2i+2k)!}{4^k\, (j+i+k)!(2i+1+k)!k!} \nonumber \\
=& \left(-\frac{1}{4}\right)^{i+j}\frac{r^{-2i}(1+2i)(2i+2j)!}{(j-i)!(j+i-1)!(2i)!}\sum_{k=0}^{j-i}\frac{(-(j-i))_{k}(j+i+1)_{k}\left(j+i+\frac{1}{2}\right)_{k}}{k!(j+i)_{k+1}(2i+1)_{k+1}}\nonumber \\
=& \left(-\frac{1}{4}\right)^{i+j}\frac{r^{-2i}(2i+2j)!}{(j-i)!(j+i)!(2i)!}\sum_{k=0}^{j-i}\frac{\left(-(j-i)\right)_{k}\left(j+i+\frac{1}{2}\right)_{k}}{k!(2i+2)_{k}} \nonumber
\end{align}
The {\it hypergeometric function} $F$ can be expressed by 
\begin{equation*}
F(-m, \beta; \gamma; z)=\sum_{k=0}^{m}\frac{(-m)_{k}(\beta)_{k}}{(\gamma)_{k}}\frac{z^{k}}{k!},
\end{equation*}
where $m$ is a positive integer and $\gamma$ is neither zero nor a negative integer (see, e.g., \cite{polyanin2001handbook}).  Therefore, 
\begin{align}
\left(\tilde{G}_{r}^{-1}\right)_{2i, 2j}
=& \left(-\frac{1}{4}\right)^{i+j}\frac{r^{-2i}(2i+2j)!}{(j-i)!(j+i)!(2i)!}F\left(-(j-i), \left(j+i+\frac{1}{2}\right) ; (2i+2); 1\right) \nonumber.
\end{align}
The special value of $F(\alpha, \beta; \gamma, z)$ at $z=1$ can be expressed in terms of a Gamma function \cite{zhang1996computation}.  That is, $F(\alpha, \beta; \gamma, 1)=\left(\Gamma(\gamma)\Gamma(\gamma-\alpha-\beta)\right) /\left(\Gamma(\gamma-\alpha)\Gamma(\gamma-\beta)\right)$.  Thus, we deduce that
\begin{align}
\left(\tilde{G}_{r}^{-1}\right)_{2i, 2j}
=& \left(-\frac{1}{4}\right)^{i+j}\frac{r^{-2i}(2i+2j)!}{(j-i)!(j+i)!(2i)!}\frac{\Gamma(2i+2)\Gamma(\frac{3}{2})}{\Gamma(i+j+2)\Gamma(i-j+\frac{3}{2})} \nonumber \\
=& \left(\left(-\frac{1}{4}\right)^{i+j}\frac{(2i+2j)!}{(j+i)!(j+i)!}\right)\left( \frac{r^{-2i}(2i+1)}{(i+j+1)}\right) \left(\frac{\Gamma(\frac{3}{2})}{(j-i)!\Gamma(i-j+\frac{3}{2})}\right). 
\end{align}

\subsection{The Odd Rows}
Similarly, from \eqref{OddrowsofInverse} in Lemma ~\ref{lem:InverseMatrix}, the nonzero entries in the odd rows of the inverse matrix $\tilde{G}_{r}^{-1}$ for $i \leq j \leq \left \lfloor \frac{N}{2} \right\rfloor$ can be rewritten as 
\begin{align}
\left(\tilde{G}_{r}^{-1}\right)_{2i+1, 2j+1} 
= & \frac{(-1)^{i+j}}{4^{i+j+1}}r^{-2i-1}(2+2i)\sum_{k=0}^{j-i}\left(-\frac{1}{4}\right)^{k} \frac{(2j+2i+2k+2)!}{\left((j-i)-k\right)!(j+i+1+k)!(2i+2+k)!k!}. \nonumber
\end{align}
By using \eqref{identities} we deduce that
\begin{align}
\left(\tilde{G}_{r}^{-1}\right)_{2i+1, 2j+1} =& \frac{(-1)^{i+j}}{4^{i+j+1}}\frac{r^{-2i-1}(2+2i)}{(j-i)!} \sum_{k=0}^{j-i}\frac{\left(-(j-i)\right)_{k}(2(j+i+1)+2k)!}{4^k\, (j+i+1+k)!(2i+2+k)!k!} \nonumber \\
=& \frac{(-1)^{i+j}}{4^{i+j+1}}\frac{r^{-2i-1}(2+2i)(2(j+i+1))!}{(j-i)!(j+i)!(2i+1)!}\sum_{k=0}^{j-i}\frac{(-(j-i))_{k}(j+i+2)_{k}\left(j+i+\frac{3}{2}\right)_{k}}{k!(j+i+1)_{k+1}(2i+2)_{k+1}}\nonumber \\
=& \frac{(-1)^{i+j}}{4^{i+j+1}}\frac{r^{-2i-1}(2(j+i+1))!}{(j-i)!(j+i+1)!(2i+1)!}\sum_{k=0}^{j-i}\frac{\left(-(j-i)\right)_{k}\left(j+i+\frac{3}{2}\right)_{k}}{k!(2i+3)_{k}} \nonumber\\
=& \frac{(-1)^{i+j}}{4^{i+j+1}}\frac{r^{-2i-1}(2(j+i+1))!}{(j-i)!(j+i+1)!(2i+1)!} F\left(-(j-i), \left(j+i+\frac{3}{2}\right) ; (2i+3); 1\right), \nonumber
\end{align}
where $F$ is a hypergeometric function. 
Again, by using a Gamma function for the special value of $F$ at $z=1$, we finally deduce that
\begin{align}
\left(\tilde{G}_{r}^{-1}\right)_{2i+1, 2j+1}
=& \frac{(-1)^{i+j}}{4^{i+j+1}}\frac{r^{-2i-1}(2(j+i+1))!}{(j-i)!(j+i+1)!(2i+1)!}\frac{\Gamma(2i+3)\Gamma(\frac{3}{2})}{\Gamma(i+j+3)\Gamma(i-j+\frac{3}{2})} \nonumber \\
=& \left(\frac{(-1)^{i+j}}{4^{i+j+1}}\frac{(2(j+i+1))!}{(j+i+1)!(j+i+1)!}\right)\left( \frac{r^{-2i-1}(2i+2)}{(i+j+2)}\right) \left(\frac{\Gamma(\frac{3}{2})}{(j-i)!\Gamma(i-j+\frac{3}{2})}\right). 
\end{align}

%
%

\end{document}